\documentclass[aos,MSNbibl,seceqn,citesort,dvips]{arximspdf}
\usepackage{graphicx}

\doi{10.1214/11-AOS906}
\volume{39}
\issue{5}
\pubyear{2011}
\firstpage{2448}
\lastpage{2476}

\makeatletter

\newtheorem{Theor}{Theorem}[section]
\newtheorem{Prop}{Proposition}[section]
\newtheorem{Lem}{Lemma}[section]
\newtheorem{Corol}{Corollary}[section]

\newproclaim{assumption}{Assumption}[section]

\newcommand{\R}{\mathbb R}
\newcommand{\N}{\mathbb N}

\newcommand{\ny}{n\to\infty}

\makeatother

\begin{document}
\begin{frontmatter}

\title{Semiparametrically efficient inference
based on signed ranks in symmetric independent component~models}
\runtitle{Rank-based inference in symmetric IC models}

\begin{aug}
\author[A]{\fnms{Pauliina} \snm{Ilmonen}\thanksref{t1}\ead[label=e1]{Pauliina.Ilmonen@gmail.com}}
\and
\author[B]{\fnms{Davy} \snm{Paindaveine}\corref{}\thanksref{t2}\ead [label=e2]{dpaindav@ulb.ac.be}
\ead[label=u1,url]{http://homepages.ulb.ac.be/\textasciitilde dpaindav}}
\runauthor{P. Ilmonen and D. Paindaveine}
\affiliation{Aalto University, University of Tampere, Universit\'{e} Libre
de Bruxelles,\break and Universit\'{e} Pierre et Marie Curie}
\address[A]{Aalto University School of Economics\\
Quantitative Methods in Economics\\
FI-00076 Aalto\\
Finland\\
\printead{e1}}
\address[B]{E.C.A.R.E.S.\\
\quad and D\'{e}partement de Math\'{e}matique\\
Universit\'{e} Libre de Bruxelles\\
50, Avenue F.D. Roosevelt, CP114/04\\
B-1050 Brussels\\
Belgium\\
\printead{e2}\\
\printead{u1}} 
\end{aug}

\thankstext{t1}{Supported by the Academy of Finland.}
\thankstext{t2}{Supported by an A.R.C. contract of the
Communaut\'{e} Fran\c{c}aise de Belgique. Davy Paindaveine is also
member of ECORE, the association between CORE and ECARES.}

\received{\smonth{1} \syear{2011}}
\revised{\smonth{6} \syear{2011}}

%
\begin{abstract}
We consider semiparametric location-scatter models for which the
\mbox{$p$-variate} observation is obtained as $X=\Lambda Z+\mu$, where $\mu$
is a $p$-vector, $\Lambda$ is a full-rank $p\times p$ matrix and the
(unobserved) random $p$-vector $Z$ has marginals that are centered
and mutually independent but are otherwise unspecified. As in blind
source separation and independent component analysis (ICA), the
parameter of interest throughout the paper is $\Lambda$. On the basis
of $n$ i.i.d. copies of $X$, we develop, under a symmetry
assumption on $Z$, \textit{signed-rank} one-sample testing and estimation
procedures for $\Lambda$. We exploit the uniform local and asymptotic
normality (ULAN) of the model to define signed-rank procedures that
are semiparametrically efficient under correctly specified densities.
Yet, as is usual in rank-based inference, the proposed procedures
remain valid (correct asymptotic size under the null, for hypothesis
testing, and root-$n$ consistency, for point estimation) under a very
broad range of densities. We derive the asymptotic properties of the
proposed procedures and investigate their finite-sample behavior
through simulations.
\end{abstract}

%

\end{frontmatter}

\section{Introduction}\label{intro}

In multivariate statistics, concepts of location and scatter are
usually defined through affine transformations of a noise vector. To be
more specific, assume that the observation $X$ is obtained through
%
%
\begin{equation}
\label{model}
X=\Lambda Z+\mu,
\end{equation}
where $\mu$ is a $p$-vector, $\Lambda$ is a full-rank $p\times p$
matrix and $Z$ is some \textit{standardized} random vector. The exact
nature of the resulting \textit{location} parameter $\mu$,
\textit{mixing matrix} parameter $\Lambda$, and \textit{scatter}
parameter $\Sigma=\Lambda\Lambda'$ crucially depends on the
standardization adopted.

The most classical assumption on $Z$ specifies that $Z$ is standard
$p$-normal. Then $\mu$ and $\Sigma$ simply coincide with the mean
vector $\mathrm{E}[X]$ and variance--covariance matrix $\operatorname{Var}[X]$
of $X$, respectively. In robust statistics, it is often rather assumed
that $Z$ is spherically symmetric about the origin of $\R^p$---in the
sense that the distribution of $OZ$ does not depend on the orthogonal
$p\times p$ matrix $O$. The resulting model in (\ref{model}) is then
called the \textit{elliptical }model. If $Z$ has finite second-order
moments, then $\mu=\mathrm{E}[X]$ and $\Sigma=c\operatorname{Var}[X]$ for
some $c>0$, but (\ref{model}) allows to define $\mu$ and $\Sigma$ in the absence of
any moment assumption.

This paper focuses on an alternative standardization of $Z$, for
which $Z$ has mutually independent marginals with common median zero.
The resulting model in (\ref{model})---the \textit{independent
component}
(\textit{IC}) \textit{model}, say---is more flexible than the elliptical model, even if
one restricts, as we will do, to vectors $Z$ with symmetrically
distributed marginals. The IC model indeed allows for heterogeneous
marginal distributions for $X$, whereas, in contrast, marginals in the
elliptical model all share---up to location and scale---the same
distribution, hence also the same tail weight. This severely affects
the relevance of elliptical models for practical applications,
particularly so for moderate to large dimensions, since it is then very
unlikely that all variables share, for example, the same tail weight.

The IC model provides the most standard setup for \textit{independent
component analysis} (\textit{ICA}), in which the mixing matrix $\Lambda$ is to
be estimated on the basis of $n$ independent copies $X_1,\ldots,X_n$
of $X$, the objective being to recover (up to a translation) the
original unobservable independent signals $Z_1,\ldots,Z_n$ by
premultiplying the $X_i$'s with the resulting $\hat\Lambda^{-1}$. It is
well known in ICA, however, that $\Lambda$ is severely unidentified:
for any $p\times p$ permutation matrix $P$ and any full-rank diagonal
matrix $D$, one can always write
%
%
\begin{equation}
\label{identif}
X
= [\Lambda P D] [(P D)^{-1} Z]+\mu
= \tilde\Lambda\tilde{Z}+\mu,
\end{equation}
where $\tilde{Z}$ still has independent marginals with median zero.
Provided that $Z$ has at most one Gaussian marginal, two matrices
$\Lambda_1$ and $\Lambda_2$ may lead to the same distribution for $X$
in (\ref{model}) if and only if they are equivalent \mbox{(we will
write $\Lambda_1\sim\Lambda_2$)} in the sense that $\Lambda
_2=\Lambda_1
PD$ for some matrices $P$ and $D$ as in (\ref{identif}); see, for
example, \cite{Th04}. In other words, under the assumption that $Z$
has at most one Gaussian marginal, permutations ($P$), sign changes and
scale transformations ($D$) of the independent components are the only
sources of unidentifiability for $\Lambda$.

This paper considers inference on the mixing matrix $\Lambda$. More
precisely, because of the identifiability issues above, we rather
consider a normalized version $L$ of $\Lambda$, where $L$ is a
well-defined representative of the class of mixing matrices that are
equivalent to $\Lambda$. This parameter $L$ is actually the parameter
of interest in ICA: an estimate of $L$ will indeed allow one to recover
the independent signals $Z_1,\ldots,Z_n$ equally well as an estimate of
any other $\Lambda$ with $\Lambda\sim L$. Interestingly, the situation
is extremely similar when considering inference on $\Sigma$ in the
elliptical model. There, $\Sigma$ is only identified up to a positive
scalar factor, and it is often enough to focus on inference about the
well-defined \textit{shape} parameter $V=\Sigma/(\det\Sigma)^{1/p}$
(e.g., in PCA, principal directions, proportions of explained
variance, etc. can be computed from $V$). Just as $L$ is a normalized
version of $\Lambda$ in the IC model, $V$ is a normalized version of
$\Sigma$ in the elliptical model, and in both classes of models, the
normalized parameters actually are the natural parameters of interest
in many inference problems. The similarities further extend to the
semiparametric nature of both models: just as the density \mbox{$g_{\|\cdot
\|
}$} of $\|Z\|$ in the elliptical model, the pdf $g_r$ of the various
independent components $Z_r$, $r=1,\ldots,p$, in the IC model, can
hardly be assumed to be known in practice.

These strong similarities motivate the approach we adopt in this paper:
we plan to conduct inference on $L$ (hypothesis testing and point
estimation) in the IC model by adopting the methodology that proved
extremely successful in \cite{Ha06c,Ha06} for inference on $V$ in the
elliptical model. This methodology combines semiparametrically
efficient inference and \textit{invariance arguments}.
In the IC model, the fixed-$(\mu,\Lambda)$ nonparametric submodels
(indexed by $g_1,\ldots,g_p$) indeed enjoy a strong invariance
structure that is parallel to the one of the corresponding elliptical
submodels (indexed by $g_{\|\cdot\|}$). As in \cite{Ha06c,Ha06}, we
exploit this invariance structure through a general result from \cite
{HW03} that allows one to derive invariant versions of efficient
central sequences, on the basis of which one can define
semiparametrically efficient (at fixed target densities $g_r=f_r$,
$r=1,\ldots,p$) invariant procedures.
As the maximal invariant associated with the invariance structure
considered turns out to be the vector of marginal signed ranks of the
residuals, the proposed procedures are of a signed-rank nature and do
not require to estimate densities.
While they achieve semiparametric efficiency under correctly specified
densities, they remain valid (correct asymptotic size under the null,
for hypothesis testing, and root-$n$ consistency, for point estimation)
under misspecified densities.

We will consider the problem of estimating $L$ and that of testing the
null $\mathcal{H}_0\dvtx L=L_0$ against the alternative $\mathcal
{H}_1\dvtx L\neq
L_0$, for some fixed $L_0$. While point estimation is undoubtedly of
primary importance for applications (e.g., in blind source separation),
one might question the practical relevance of the testing problem
considered, especially when $L_0$ is not the $p$-dimensional identity
matrix. 
Solving this generic testing problem, however, is the main step in
developing tests for any linear hypothesis on $L$, and we will
explicitly describe the resulting tests in the sequel. An extensive
study of these tests is beyond the scope of the present paper, though;
we refer to \cite{IC} for an extension of our tests to the particular
case of testing the (linear) hypothesis that $L$ is block-diagonal, a
problem that is obviously important in practice (nonrejection of the
null would indeed allow practitioners to proceed with two separate,
lower-dimensional, analyses). Testing linear hypotheses on $L$ includes
many other testing problems of high practical relevance, such as
testing that a given column of $L$ is equal to some fixed $p$-vector,
and testing that a given entry of $L$ is zero---the practical
importance of these two testing problems, in relation, for example,
with functional magnetic resonance imaging (fMRI), is discussed
in \cite{Ol11}.

The paper is organized as follows. In Section \ref{modelULANsec}, we
fix the notation and describe the model (Section \ref{modelsec}), state
the corresponding \textit{uniformly locally and asymptotically normal}
(\textit{ULAN}) property that allows us to determine semiparametric
efficiency bounds (Section \ref{ULANsec}) and then introduce, in
relation with invariance arguments, \textit{rank-based} efficient
central sequences (Section \ref{invarianceranksec}). In Sections
\ref{testsec} and \ref {estimating}, we develop the resulting rank
tests and estimators for the mixing matrix $L$, respectively. Our
estimators actually require the delicate estimation of $2p(p-1)$
``cross-information coefficients,'' an issue we solve in Section
\ref{estimcrosssec} by generalizing the method recently developed in
\cite{r12}. In Section \ref{simusec}, simulations are conducted both to
compare the proposed estimators with some competitors and to
investigate the validity of asymptotic results---simulation results for
hypothesis testing are provided in the supplementary article
\cite{IP11a}. Finally, the \hyperref[app]{Appendix} states some
technical results (Appendix \ref{appensec}) and reports proofs
(Appendix \ref{appensecB}).


\section{The model, the ULAN property and invariance arguments}
\label{modelULANsec}

\subsection{The model}\label{modelsec}

As we already explained, the IC model above
suffers from severe identifiability issues for $\Lambda$. To solve
this, we map each $\Lambda$ onto a\vspace*{1pt} unique representative
$L=\Pi (\Lambda )$ of the collection of mixing matrices $\tilde\Lambda$
that satisfy $\tilde\Lambda\sim\Lambda$ (the equivalence class of
$\Lambda$ for $\sim$).
We propose the mapping
\[
\Lambda\mapsto\Pi(\Lambda)=\Lambda D^+_{1} P D_2,
\]
where $D^+_{1}$ is the positive definite diagonal matrix that makes
each column of $\Lambda D^+_{1}$ have Euclidean norm one, $P$ is the
permutation matrix for which the matrix $B=(b_{ij})=\Lambda D^+_{1} P$
satisfies $|b_{ii}|>|b_{ij}|$ for all $i<j$
and $D_2$ is the diagonal matrix such that all diagonal entries of $\Pi
(\Lambda)=\Lambda D^+_{1} P D_2$ are equal to one.

If one restricts to the collection $\mathcal{M}_p$ of mixing
matrices $\Lambda$ for which no ties occur in the permutation step
above, it can easily be shown that, for any $\Lambda_1,\Lambda_2\in
\mathcal{M}_p$, we have that $\Lambda_1\sim\Lambda_2$ iff $\Pi
(\Lambda
_1)=\Pi(\Lambda_2)$, so that this mechanism succeeds in identifying a
unique representative in each class of equivalence (this is ensured
with the double scaling scheme above, which may seem a bit complicated
at first). Besides, $\Pi$ is then a continuously differentiable mapping
from $\mathcal{M}_p$ onto $\mathcal{M}_{1p}:=\Pi(\mathcal{M}_{p})$.
While ties may always be taken care of in some way (e.g., by basing the
ordering on subsequent rows of the matrix $B$), they may prevent the
mapping $\Pi$ to be continuous, which would cause severe problems and
would prevent us from using the Delta method in the sequel. It is
clear, however, that the restriction to $\mathcal{M}_p$ only gets rid
of a few particular mixing matrices, and will not have any implications
in practice.

The parametrization of the IC model we consider is then associated with
%
%
\begin{equation}
\label{model2}
X=L Z+\mu,
\end{equation}
where $\mu\in\R^p$, $L\in\mathcal{M}_{1p}$ and $Z$ has independent
marginals with common median zero. Throughout, we further assume that
$Z$ admits a density with respect to the Lebesgue measure on $\R^p$,
and that it has $p$ symmetrically distributed marginals, among which at
most one is Gaussian (as explained in the \hyperref
[intro]{Introduction}, this limitation
on the number of Gaussian components is needed for $L$ to be
identifiable). We will denote by $\mathcal{F}$ the resulting collection
of densities for $Z$. Of course, any $g\in\mathcal{F}$ naturally
factorizes into $g(z)=\prod_{r=1}^p g_r(z_r)$, where $g_r$ is the
symmetric density of~$Z_r$.

%

The hypothesis under which $n$ mutually independent observations $X_i$,
$i=1,\ldots,n$, are obtained from (\ref{model2}), where $Z$ has
density $g\in\mathcal{F}$, will be denoted
as $\mathrm{P}^{(n)}_{\vartheta,g}$, with $\vartheta=(\mu',(
\operatorname{vecd}^\circ
L)')'\in\Theta=\R^p\times\operatorname{vecd}^\circ(\mathcal{M}_{1p})$,
or alternatively, as $\mathrm{P}^{(n)}_{\mu,L,g}$; for any\vspace*{1pt}
$p\times p$ matrix $A$, we write $\operatorname{vecd}^\circ A$ for the
$p(p-1)$-vector obtained by removing the $p$ diagonal entries of $A$
from its usual vectorized form $\operatorname{vec} A$ (diagonal entries
of $L$ are all equal to one, hence should not be included in the
parameter).

The resulting semiparametric model is then
%
%
\begin{equation} \label{semimodel}
\mathcal{P}^{(n)} := \bigcup_{g\in\mathcal{F}} \mathcal{P}_g^{(n)} :=
\bigcup_{g\in\mathcal{F}} \bigcup_{\vartheta\in\Theta} \bigl\{
\mathrm{P}^{(n)}_{\vartheta,g}\bigr\}.
\end{equation}
Performing semiparametrically efficient inference on $\vartheta$, at a
fixed $f\in\mathcal{F}$, typically requires that the corresponding
parametric submodel $\mathcal{P}_f^{(n)}$ satisfies the \textit{uniformly
locally and asymptotically normal} (\textit{ULAN}) property.

\subsection{The ULAN property}
\label{ULANsec}

As always, the ULAN property requires technical regularity conditions
on $f$. In the present context, we need that each corresponding
univariate pdf $f_r$, $r=1,\ldots,p$, is absolutely continuous (with
derivative $f'_r$, say) and satisfies
\[
\sigma^2_{f_r}:=\int_{-\infty}^\infty y^2 f_r(y) \,dy < \infty,\qquad
\mathcal{I}_{f_r}:=\int_{-\infty}^\infty\varphi_{f_r}^2(y)
f_r(y) \,dy
< \infty
\]
and
\[
\mathcal{J}_{f_r}:=\int_{-\infty}^\infty y^2\varphi_{f_r}^2(y)
f_r(y)
\,dy < \infty,
\]
where we let $\varphi_{f_r}:=-f'_r/f_r$.
In the sequel, we denote by $\mathcal{F}_{\mathrm{ulan}}$ the
collection of
pdfs $f\in\mathcal{F}$ meeting these conditions.

For any $f\in\mathcal{F}_{\mathrm{ulan}}$, let
$\gamma_{rs}(f):=\mathcal{I}_{f_r}\sigma^2_{f_s}$, define the optimal
$p$-variate location score function $\varphi_{f}\dvtx\R^p\to\R^p$
through $z=(z_1,\ldots,z_p)'\mapsto\varphi_{f}
(z)=(\varphi_{f_1}(z_1),\break \ldots, \varphi_{f_p}(z_p))'$, and denote by $
\mathcal{I}_{f}$
the diagonal matrix with diagonal entries $\mathcal{I}_{f_r}$,
$r=1,\ldots,p$.
Further write $I_\ell$ for the $\ell$-dimensional identity matrix and define
%
\vspace*{-6pt}
\[
C :=
\sum_{r=1}^p
\sum_{s=1}^{p-1}
\bigl(e_r e_r'\otimes u_s e_{s+\delta_{[s\geq r]}}'\bigr),
\]
where $\otimes$ is the usual Kronecker product, $e_r$ and $u_r$ stand
for the $r$th vectors of the canonical basis of $\mathbb R^p$
and $\mathbb R^{p-1}$, respectively, and $\delta_{[s\geq r]}$ is equal
to one if $s\geq r$ and to zero otherwise. The following ULAN result
then easily follows from Proposition 2.1 in \cite{IC} by using a simple
chain rule argument.

\vspace*{-3pt}
\begin{Prop}\label{prulan}
Fix $f\in\mathcal{F}_{\mathrm{ulan}}$. Then
the collection of probability distributions
$
\mathcal{P}_f^{(n)}
$
is ULAN, with central sequence
%
%
\begin{equation}
\label{centralse}\qquad
\Delta_{\vartheta,f}
=
\pmatrix{
\Delta_{\vartheta,f;1}\cr
\Delta_{\vartheta,f;2}}
=
\pmatrix{
\displaystyle n^{-1/2}(L^{-1})'\sum_{i=1}^n\varphi_{f}(Z_i)\vspace*{2pt}\cr
\displaystyle n^{-1/2}C(I_p\otimes L^{-1})'\sum_{i=1}^n\operatorname{vec}\bigl(\varphi
_{f}(Z_i)Z_i'-I_p\bigr)},
\end{equation}
where $Z_i=Z_i(\vartheta)=L^{-1}(X_i-\mu)$, and full-rank information matrix
\[
\Gamma_{L,f}
=
\pmatrix{
\Gamma_{L,f;1} & 0 \cr
0 & \Gamma_{L,f;2}},
\]
where
$
\Gamma_{L,f;1}
:= (L^{-1})'\mathcal{I}_{f}L^{-1}
$
and
\begin{eqnarray*}
\Gamma_{L,f;2}
&:=&
C (I_p\otimes L^{-1})'
\Biggl[
\sum_{r=1}^p
(\mathcal{J}_{f_r}-1)(e_re_r'\otimes e_re_r')
\\[-3pt]
& &\hphantom{C (I_p\otimes L^{-1})'
\Biggl[}
{}
+ \sum_{r,s=1, r\neq s}^p
\bigl(
\gamma_{sr}(f)(e_re_r'\otimes e_s e_s') + (e_re_s'\otimes e_s e_r')
\bigr)
\Biggr]\\[-3pt]
&&{}\times
(I_p\otimes L^{-1})C'.
\end{eqnarray*}
More precisely, for any $\vartheta_n=\vartheta+ O(n^{-1/2})$ $($with
$\vartheta=(\mu',(\operatorname{vecd}^\circ L)')')$ and any bounded sequence
$(\tau_n)$ in $\mathbb R^{p^2}$, we have that, under $\mathrm{P}^{(n)}
_{\vartheta_n,f}$ as $n \rightarrow\infty$,
\[
\log\bigl(d\mathrm{P}^{(n)}_{\vartheta_n + n^{-1/2}\tau_n,f} /d\mathrm
{P}^{(n)}_{\vartheta
_n,f}\bigr)=\tau_n'\Delta_{\vartheta_n,f}-\tfrac{1}{2}\tau_n'\Gamma
_{L,f}\tau
_n + o_{\mathrm{P}}(1),
\]
and $\Delta_{\vartheta_n,f}$ converges in distribution to a
$p^2$-variate normal distribution with mean zero and covariance matrix
$\Gamma_{L,f}$.
\end{Prop}


Semiparametrically efficient (at $f$) inference procedures on $L$ then
may be based on the so-called
\textit{efficient central sequence} $\Delta^{*}_{\vartheta,f;2}$
resulting from $\Delta_{\vartheta,f;2}$ by performing adequate
tangent space projections; see \cite{Bi93}. Under $\mathrm
{P}^{(n)}_{\vartheta
,f}$, $\Delta^{*}_{\vartheta,f;2}$ is still asymptotically normal with
mean zero, but now with covariance matrix $\Gamma^*_{L,f;2}$ \label
{gammatoto} (the \textit{efficient information matrix}). This
matrix $\Gamma^*_{L,f;2}$ settles the semiparametric efficiency bound
at $f$ when performing inference on $L$. For instance, an
estimator $\hat{L}$ is semiparametrically efficient at $f$ if
%
%
\begin{equation}\label{estimsemiparam}
\sqrt{n} \operatorname{vecd}^\circ(\hat{L}-L)
\stackrel{\mathcal{L}}{\to} \mathcal{N}_{p(p-1)} ( 0,
(\Gamma^*_{L,f;2})^{-1} ).
\end{equation}
The performance of semiparametrically efficient tests on $L$ can
similarly be characterized in terms of $\Gamma^*_{L,f;2}$: 
a test of $\mathcal{H}_0\dvtx L=L_0$ is semiparametrically efficient at $f$
(at asymptotic level $\alpha$) if its asymptotic powers under local
alternatives of the form $\mathcal{H}_1^{(n)}\dvtx L=L_0+n^{-1/2}H$,
where $H$
is an arbitrary $p\times p$ matrix with zero diagonal entries, are
given by
%
%
\begin{equation}
\label{noncentrsemiparam}
1 - \Psi_{p(p-1)} \bigl( \chi^2_{p(p-1),1-\alpha}
; (\operatorname{vecd}^\circ H)' \Gamma^*_{L_0,f;2}
(\operatorname{vecd}^\circ H) \bigr),
\end{equation}
where $\chi^2_{p(p-1),1-\alpha}$ stands for the $\alpha$-upper quantile
of the $\chi^2_{p(p-1)}$ distribution, and $\Psi_{p(p-1)}( \cdot
;\delta)$ denotes the cumulative distribution function of the
noncentral $\chi^2_{p(p-1)}$ distribution with \label
{pagechideuxnoncentr} noncentrality parameter $\delta$.

\subsection{Invariance arguments}
\label{invarianceranksec}

Instead of the classical tangent space projection approach to
compute $\Delta^{*}_{\vartheta,f;2}$ (as in \cite{CB06}), we adopt an
approach---due to \cite{HW03}---that rather exploits the invariance
structure of the model considered. This will provide a version of the
efficient central sequence (parallel to central sequences, efficient
central sequences are defined up to $o_\mathrm{P}(1)$'s only) that is based
on \textit{signed ranks}. Here, signed ranks are defined as
$S_i(\vartheta)=(S_{i1}(\vartheta),\ldots,S_{ip}(\vartheta))'$ and
$R^+_i(\vartheta)=(R^+_{i1}(\vartheta),\ldots,R^+_{ip}(\vartheta))'$,
where $S_{ir}(\vartheta)$ is the sign of $Z_{ir}(\vartheta)=(L^{-1}
(X_i-\mu))_{r}$ and $R^+_{ir}(\vartheta)$ is the rank of
$|Z_{ir}(\vartheta)|$ among
$|Z_{1r}(\vartheta)|,\ldots,|Z_{nr}(\vartheta)|$. This signed-rank
efficient central sequence---$\underline{\Delta}^{*}_{\vartheta,f;2}$,
say---is given in Theorem \ref{th1} below (the asymptotic behavior
of $\underline{\Delta}^{*}_{\vartheta,f;2}$ will be studied in
Appendix \ref{appensec}).\vspace*{1pt}

To be able to state Theorem \ref{th1}, we need to introduce the
following notation.
Let
$z\mapsto F_+
(z)=(F_{+1}(z_1), \ldots, F_{+r}(z_p))'$, with $F_{+r}(t):=\mathrm{P}^{(n)}
_{\vartheta,f}[|Z_r(\vartheta)|<t]=2(\int_{-\infty}^tf_r(s)
\,ds)-1$, $t\geq0$. Based on this,
define
$
\underline{\Delta}^{*}_{\vartheta,f;2}
:=
C(I_p\otimes L^{-1})'\*
\operatorname{vec} \underline{T}_{\vartheta,f},
$
with
\begin{eqnarray*}
\underline{T}_{\vartheta,f}
&:=&
\operatorname{odiag} \Biggl[ \frac{1}{\sqrt{n}}
\sum_{i=1}^n
\biggl(S_i(\vartheta)\odot\varphi_{f}\biggl(F_+^{-1}\biggl(\frac
{R^+_{i}(\vartheta)}{n+1}\biggr)\biggr)\biggr)\\
&&\hphantom{\operatorname{odiag} \Biggl[ \frac{1}{\sqrt{n}}
\sum_{i=1}^n}\hspace*{10.7pt}
{}\times
\biggl(S_i(\vartheta)\odot F_+^{-1}\biggl(\frac{R^+_{i}(\vartheta
)}{n+1}\biggr)\biggr)'
\Biggr],
\end{eqnarray*}
where $\odot$ is the Hadamard (i.e., entrywise) product of two vectors,
and where $\operatorname{odiag}(A)$ denotes the matrix obtained from
$A$ by
replacing all diagonal entries with zeros. Finally, let $\underline
{\mathcal{F}}_{\mathrm{ulan}}$ be the collection of pdfs $f\in\mathcal
{F}_{\mathrm{ulan}}$ for which each $\varphi_{f_r}$, $r=1,\ldots,p$, is
continuous and can be written as the difference of two monotone
increasing functions. We then have the following result (see
Appendix \ref{appensecB} for a proof).
\begin{Theor}\label{th1}
Fix $\vartheta=(\mu',(\operatorname{vecd}^\circ L)')'\in\Theta$ and
$f\in
\underline{\mathcal{F}}_{\mathrm{ulan}}$.
Then, \textup{(i)} denoting by $\mathrm{E}^{(n)}_{\vartheta,f}$
expectation under $\mathrm{P}^{(n)}_{\vartheta,f}$,
\begin{eqnarray*}
\underline{\Delta}^{*}_{\vartheta,f;2}
:\!&=&
C(I_p\otimes L^{-1})'
\operatorname{vec} \underline{T}_{\vartheta,f}
\\
&=&
\mathrm{E}^{(n)}_{\vartheta,f}[\Delta_{\vartheta,f;2} |
S_1(\vartheta),
\ldots, S_n(\vartheta), R^+_1(\vartheta), \ldots, R^+_n(\vartheta)]
+ o_{L^2}(1)
\end{eqnarray*}
as $\ny$, under $\mathrm{P}^{(n)}_{\vartheta,f}$; \textup{(ii)} the signed-rank
quantity $\underline{\Delta}^{*}_{\vartheta,f;2}$ is a version of
the efficient
central sequence at $f$ [i.e., $\underline{\Delta }^{*}_{\vartheta
,f;2}=\Delta^{*}_{\vartheta,f;2}+o_{L^2}(1)$ as $\ny$, under
$\mathrm{P}^{(n)} _{\vartheta,f}$].
\end{Theor}

Would
the (nonparametric) fixed-$\vartheta$ submodels
$
\mathcal{P}_\vartheta^{(n)}
:=
\bigcup_{g\in\mathcal{F}}
\{ \mathrm{P}^{(n)}_{\vartheta,g}\}
$
of the semiparametric
model $\bigcup_{\theta\in\Theta} \bigcup_{g\in\mathcal{F}}
\{\mathrm{P}^{(n)}_{\theta,g}\}$ in (\ref{semimodel}) be invariant
under a group of transformations $ \mathcal{G}^\vartheta $ that
generates $\mathcal{P}_\vartheta^{(n)}$, then the main result of \cite
{HW03} would show that the expectation
of the original central sequence $\Delta_{\vartheta,f;2}$ conditional
upon the corresponding maximal invariant---$\mathcal{I}^{(n)}_{
\max}(\vartheta)$, say---is a version of the efficient central
sequence $\Delta^{*}_{\vartheta,f;2}$ at $f$:
as $\ny$, under $\mathrm{P}^{(n)}_{\vartheta,f}$,
%
%
\begin{equation}\label{HHWW03}
\Delta^{*}_{\vartheta,f;2}
=
\mathrm{E}^{(n)}_{\vartheta,f}\bigl[\Delta_{\vartheta,f;2} | \mathcal
{I}^{(n)}_{\max}(\vartheta)\bigr]
+o_{L^2}(1).
\end{equation}

Such an invariance structure actually exists and the relevant group
$
\mathcal{G}^\vartheta$ collects all transformations
\begin{eqnarray*}
g^\vartheta_h
\dvtx
\R^p\times\cdots\times\R^p
& \to&
\R^p\times\cdots\times\R^p,
\\
(x_1,\ldots,x_n)
&\mapsto&
\bigl(L h(z_1(\vartheta))+\mu,\ldots,L h(z_n(\vartheta))+\mu\bigr),
\end{eqnarray*}
with $z_{i}(\vartheta):=L^{-1} (x_i-\mu)$
and 
$h((z_1,\ldots,z_p)^\prime)=(h_1(z_1),\ldots,h_p(z_p))^\prime$, where
each $h_r$, $r=1,\ldots,p$, is continuous, odd, monotone increasing and
fixes $+\infty$. It is easy to check that $\mathcal{P}_\vartheta
^{(n)}$ is
invariant under (and is generated by) $
\mathcal{G}^\vartheta$, and that the corresponding maximal invariant is
the vector of signed ranks
%
%
\begin{equation} \label{signedranks}
\mathcal{I}^{(n)}_{\max}(\vartheta)=(S_1(\vartheta),\ldots
,S_n(\vartheta
),R^{+}_1(\vartheta),\ldots,R^{+}_n(\vartheta));
\end{equation}
Theorem \ref{th1}(ii) then follows from (\ref{HHWW03}) and
Theorem \ref{th1}(i).

%
%

Inference procedures based on $\underline{\Delta}^{*}_{\vartheta,f;2}$,
unlike those (from \cite{CB06}) based on the efficient central
sequence\vspace*{1pt} $\Delta^{*}_{\vartheta,f;2}$ obtained through tangent space
projections, are measurable with respect to signed ranks, hence enjoy
all nice properties usually associated with rank methods: robustness,
ease of computation, validity without density estimation (and, for
hypothesis testing, even distribution-freeness), etc.



\section{Hypothesis testing}\label{testsec}

We now consider the problem of testing the null hypothesis $\mathcal
{H}_0\dvtx L=L_0$ against the alternative $\mathcal{H}_1\dvtx L\neq
L_0$, with
unspecified underlying density $g$. Beyond their intrinsic interest,
the resulting tests will play an important role in the construction of
the $R$-estimators of Section \ref{estimating} below, and they pave the
way to testing linear hypotheses on $L$.

The objective here is to define a test that is semiparametrically
efficient at some target density $f$, yet that remains valid---in the
sense that it meets asymptotically the level constraint---under a very
broad class of densities $g$. As we will show, this objective is
achieved\vadjust{\goodbreak} by the signed-rank test---$\underline{\phi}_f$, say---that
rejects $\mathcal{H}_0$ at asymptotic level $\alpha\in(0,1)$ whenever
%
%
\begin{equation} \label{teststatorig}
\underline{Q}_f := (\underline{\Delta}^{*}_{\hat\vartheta_{0},f;2})'
(\Gamma
^*_{L_0,f;2})^{-1}\underline{\Delta}^{*}_{\hat\vartheta_{0},f;2}
>
\chi^2_{p(p-1),1-\alpha},
\end{equation}
where $\Gamma^*_{L,f;2}$ was introduced on Page \pageref{gammatoto}
(an explicit expression is given below)
and where $\hat\vartheta_{0}=(\hat\mu',(\operatorname{vecd}^\circ L_0)')'$
is based on a sequence of estimators $\hat\mu$ that is locally
asymptotically discrete (see Appendix \ref{appensec} for a precise
definition) and root-$n$ consistent under the null.

Possible choices for $\hat\mu$ include (discretized versions of) the
sample mean $\bar{X}:=\frac{1}{n}\sum_{i=1}^n X_i$ or the
transformation-retransformation componentwise median \label
{componentmed} $\hat\mu_{\mathrm{Med}}:=L_0 \operatorname
{Med}[L_0^{-1}X_1,\ldots
,L_0^{-1}X_n]$, where $\operatorname{Med}[\cdot]$ returns the vector of
univariate medians. We favor the sign estimator $\hat\mu_{\mathrm{Med}}$,
since it is very much in line with the signed-rank tests $\underline
{\phi}_f$ and enjoys good robustness properties.
However, we stress that Theorem \ref{theotest} below, which states the
asymptotic properties of the proposed signed-rank tests, implies that
the choice of $\hat\mu$ does not affect the asymptotic properties
of $\underline{\phi}_f$, at any $g\in{\mathcal{F}}_{\mathrm{ulan}}$.

In order to state this theorem, we need to define
%
%
\begin{eqnarray}\label{definGstar}\quad
\Gamma^*_{L,f,g;2} &:=& C (I_p\otimes L^{-1})' G_{f,g} (I_p\otimes
L^{-1})C'
\nonumber\\
&:=& C (I_p\otimes L^{-1})'
\nonumber\\[-8pt]\\[-8pt]
&&{} \times\Biggl[ \sum_{r,s=1, r\neq s}^p \bigl( \gamma_{sr}(f,g)(e_re_r'\otimes
e_s e_s') + \rho_{rs}(f,g)(e_re_s'\otimes e_s e_r') \bigr) \Biggr]\nonumber\\
&&{}\times
(I_p\otimes L^{-1})C',
\nonumber
\end{eqnarray}
where we let
%
%
\begin{equation}
\label{crossinfocoeffgam}\quad
\gamma_{rs}(f,g) := \int_{0}^1
\varphi_{f_r}(F_r^{-1}(u)) \varphi_{g_r}(G_r^{-1}(u)) \,du \times
\int_{0}^1 F_s^{-1}(u) G_s^{-1}(u) \,du
\end{equation}
and
%
%
\begin{equation}
\label{crossinfocoeffrho}\qquad \rho_{rs}(f,g) := \int_{0}^1 F_r^{-1}(u)
\varphi_{g_r}(G_r^{-1}(u)) \,du \times \int_{0}^1
\varphi_{f_s}(F_s^{-1}(u)) G_s^{-1}(u) \,du.
\end{equation}
We also\vspace*{2pt} let $\Gamma^*_{L,f;2}:=\Gamma^*_{L,f,f;2}$ and
$G_{f}:=G_{f,f}$, that involve $\gamma_{rs}(f,f)=\gamma_{rs}(f)$ (see
Section \ref{ULANsec}) and $\rho_{rs}(f,f)=1$. We then have the
following result (see Appendix \ref{appensecB} for a proof).
\begin{Theor}\label{theotest}
Fix
$f\in\underline{\mathcal{F}}_{\mathrm{ulan}}$. Then,
\textup{(i)} under $\mathrm{P}^{(n)}_{\vartheta_0,g}$ and under
$\mathrm {P}^{(n)}_{\vartheta _0+n^{-1/2}\tau,g}$, with
$\vartheta_0=(\mu',(\operatorname{vecd}^\circ L_0)')'$,
$\tau=(\tau_1',\tau_2')' \in\R^p\times\R^{p(p-1)}$ and $g\in
{\mathcal{F}}_{\mathrm{ulan}}$,
%
\[
\underline{Q}_f \stackrel{\mathcal{L}}{\to} \chi^2_{p(p-1)}
\quad\mbox{and}\quad \underline{Q}_f \stackrel{\mathcal{L}}{\to} \chi^2_{p(p-1)} (\tau_2'
(\Gamma^*_{L_0,f,g;2})' (\Gamma^*_{L_0,f;2})^{-1} \Gamma^*_{L_0,f,g;2}
\tau_2),
\]
respectively, as $\ny$.
\textup{(ii)} The sequence
of tests $\underline{\phi}_f^{(n)}$ has asymptotic level $\alpha$
under $\bigcup_{\mu\in\R^p} \bigcup_{g\in{\mathcal{F}}_{\mathrm{ulan}}}
\{
\mathrm{P}^{(n)}_{\mu,L_0,g}\}$.
\textup{(iii)} The\vadjust{\goodbreak} sequence of tests
$\underline{\phi}_f^{(n)}$ is semiparamet\-ri\-cally efficient, still at
asymptotic level $\alpha$, when testing $\mathcal{H}_0\dvtx L=L_0$
against $H_1^f\dvtx L\neq L_0$ with noise density $f$ (i.e., when
testing $\bigcup_{\mu\in\R^p} \bigcup_{g\in {\mathcal
{F}}_{\mathrm{ulan}}} \{\mathrm{P}^{(n)}_{\mu,L_0,g}\}$ against
$\bigcup _{\mu \in\R^p} \bigcup_{L\in\mathcal{M}_{1p}\setminus\{L_0\}}
\{\mathrm{P}^{(n)}_{\mu ,L,f}\})$.
\end{Theor}

The test
$\underline{\phi}_f$ achieves semiparametric efficiency at $f$
[Theorem \ref{theotest}(iii)], and also at any ${f_\sigma}$,
with $f_\sigma(z):=\prod_{r=1}^p \sigma_r^{-1} f_{r}(z_r/\sigma_r)$,
where $\sigma_r>0$ for all $r$---it can indeed be checked that
$\underline{\phi}_{f_\sigma}=\underline{\phi}_{f}$.
Most importantly, Theorem \ref{theotest} shows also that $\underline
{\phi}_{f}$ remains
valid under any $g\in{\mathcal{F}}_{\mathrm{ulan}}$. By proceeding
as in Lemma 4.2 of \cite{IC}, this can even be extended to any $g\in
\mathcal{F}$, which allows us to avoid any finite moment condition.

This is to be compared to the semiparametric approach of Chen and
Bickel \cite{CB06}---these authors focus on point estimation, but their
methodology also leads to tests that enjoy the same properties as their
estimators. Their procedures achieve uniform (in $g$) semiparametric
efficiency, while our methods achieve semiparametric efficiency at the
target density $f$ only---more precisely, at any
corresponding $f_\sigma$.
However, it turns out that the performances of our procedures do not
depend much on the target density $f$, so that our procedures are close
to achieving uniform (in $g$) semiparametric efficiency; see the
simulations in the supplemental article \cite{IP11a}. As any uniformly
semiparametrically efficient procedures (see \cite{Am02}), Chen and
Bickel's procedures require estimating $g$, hence choosing various
smoothing parameters. In contrast, our procedures, by construction, are
invariant (here, signed-rank) ones. As such, 
they do not require us to estimate densities, and they are
robust, easy to compute, etc.

One might still object that the choice of $f$ is quite arbitrary. 
%
This choice should be based on the practitioner's prior belief on the
underlying densities. If he/she has no such prior belief, a kernel
estimate $\hat{f}$ of $f$ could be used. The resulting
test $\underline
{\phi}_{\hat{f}}$ would then enjoy the same properties as
any $\underline{\phi}_f$ in terms of validity, since kernel density
estimators, in the symmetric case considered, typically are measurable with respect to the order
statistics of the $|Z_{ir}(\hat\vartheta_{0})|$'s, that,\vspace*{1pt}
asymptotically, are stochastically\vspace*{1pt} independent of the signed
ranks $S_{ir}(\hat\vartheta_{0}),R^+_{ir}(\hat\vartheta_{0})$ used
in $\underline{\phi}_f$; see \cite{HW03} for details. The
test $\underline{\phi}_{\hat{f}}$ would further achieve uniform
semiparametric efficiency.

Further results on the proposed tests are given in the supplemental
article \cite{IP11a}. More precisely, a simple explicit expression of
the test statistics, local asymptotic powers of the corresponding
tests, and simulation results can be found there.

We finish this section by describing the extension of our signed-rank
tests to the problem of testing a fixed (arbitrary) linear hypothesis
on $L$, which includes many instances of high practical relevance (we
mentioned a few in the \hyperref[intro]{Introduction}).
Denoting by $\mathcal{V}(\Omega)$ the vector space that is spanned by
the columns of the $p(p-1) \times\ell$ matrix $\Omega$ (which is
assumed to have full rank $\ell$), we consider the testing
problem\looseness=-1
%
%
\begin{equation}\label{lintestprob}
\cases{
\mathcal{H}_0(L_0,\Omega)\dvtx(\operatorname{vecd}^\circ L)\in(
\operatorname{vecd}^\circ L_0)+\mathcal{V}(\Omega)
\cr
\mathcal{H}_1(L_0,\Omega)\dvtx(\operatorname{vecd}^\circ L)\notin(
\operatorname{vecd}^\circ L_0)+\mathcal{V}(\Omega)
,}
\end{equation}
for some fixed $L_0\in\mathcal{M}_{1p}$. If one forgets about the
tacitly assumed constraint that $L\in\mathcal{M}_{1p}$ in (\ref
{lintestprob}), the null hypothesis above imposes a set of linear
constraints on $L$. This clearly includes all testing problems
mentioned in the \hyperref[intro]{Introduction}: testing that a given
column of $L$ is
equal to a fixed vector, testing that a given (off-diagonal) entry
of $L$ is zero and testing block-diagonality of $L$.

Inspired by the tests from \cite{Le86} (Section 10.9), the analog of
our signed-rank test $\underline{\phi}_f$ above then
rejects $\mathcal
{H}_0(L_0,\Omega)$ for large values of
%
\[
\underline{Q}_f(L_0,\Omega) :=
(\underline{\Delta}^{*}_{\hat\vartheta,f;2})' P_\Omega
\underline{\Delta}^{*}_{\hat\vartheta,f;2}
\]
with\vspace*{-1pt}
$ P_\Omega:=(\Gamma^*_{\hat L,f;2})^{-} -
\Omega(\Omega^\prime\Gamma^*_{\hat L,f;2} \Omega)^{-} \Omega ^\prime, $
where
$B^-$ denotes the Moore--Penrose pseudoinverse of $B$, and where $\hat
\vartheta=(\hat\mu',(\operatorname{vecd}^\circ\hat{L})')'$ is an estimator
of $\vartheta$ that is locally and asymptotically discrete, root-$n$
consistent under the null, and \textit{constrained}---in the sense
that $\hat L$ satisfies the linear constraints in $\mathcal
{H}_0(L_0,\Omega)$.

It can be shown that this signed-rank test achieves semiparametric
optimality at $f$ (the relevant optimality concept here is \textit{most
stringency}; see, e.g., \cite{IC} for a discussion) and remains valid
under any $g\in\underline{\mathcal{F}}_{\mathrm{ulan}}$. Its null asymptotic
distribution is still chi-square, now with $r:=\operatorname
{Trace}[P_\Omega
\Gamma^*_{L,f;2}]$ degrees of\vspace*{1pt} freedom (this directly follows from
Theorem 9.2.1 in \cite{Ra71} and Theorem \ref{th2}); at asymptotic
level~$\alpha$, the resulting asymptotic critical value (that actually
does not depend on the true value $L$) therefore is $\chi
^2_{r;1-\alpha
}$. Just as for the tests $\underline{\phi}_f$, it is still
possible to compute asymptotic powers under sequences of local
alternatives. It is clear, however, that a thorough study of the
properties of the tests above, for a general linear hypothesis, is
beyond the scope of the present paper, hence is left for future
research. In the important particular case of testing block-diagonality
of $L$, a complete investigation of the signed-rank tests can be found
in \cite{IC}.






\section{Point estimation}\label{estimating}

We turn to the problem of estimating $L$, which is of primary
importance for applications.
Denoting by $\underline{Q}_f=\underline{Q}_f(L_0)$ the signed-rank test
statistic for $\mathcal{H}_0\dvtx L=L_0$ in (\ref{teststatorig}), a natural
signed-rank estimator of $L$ is obtained by ``inverting the
corresponding test,''
\[
\hat{\underline{L}}_{f;\arg\min}=\mathop{\arg\min}_{L\in\mathcal
{M}_{1p}}\underline{Q}_f(L).
\]
This estimator, however, is not satisfactory: as any signed-rank
quantity, the objective function $L\mapsto\underline{Q}_f(L)$ is
piecewise\vspace*{-1pt} constant, hence discontinuous and nonconvex,
which makes it very difficult to derive the asymptotic properties of
$\hat{\underline {L}}_{f;\arg\min}$. It is also virtually impossible to
compute $\hat {\underline{L}}_{f;\arg\min}$ in practice, since this
lack of smoothness and convexity essentially forces computing the
estimator by simply running over a grid of possible values of the
$p(p-1)$-dimensional parameter $L$---a strategy that cannot provide a
reasonable approximation of $\hat{\underline{L}}_{f;\arg\min}$,
even for moderate values of $p$. Finally, there is no way to estimate
the asymptotic covariance matrix of $\hat{\underline{L}}_{f;
\arg\min}$, which rules out the possibility to derive confidence zones
for $L$, hence drastically restricts the practical relevance of this
estimator.

In order to avoid the aforementioned drawbacks, we propose adopting a
one-step approach that was first used in \cite{Ha06c} for the problem
of estimating the shape of an elliptical distribution or in \cite
{Ha08b} in a more general context. The resulting one-step signed-rank
estimators---in the sequel, we simply speak of \textit{one-step rank
estimators} or \textit{one-step $R$-estimators}---can easily be computed in
practice, their asymptotic properties can be derived explicitly, and
their asymptotic covariance matrix can be estimated consistently.

\subsection{One-step $R$-estimators of $L$}
\label{estimasymptsec}

To initiate the one-step procedure, a preliminary estimator is needed.
In the present context, we will assume that a root-$n$ consistent and
locally asymptotically discrete estimator $\tilde\vartheta=(\tilde
\mu
',(\operatorname{vecd}^\circ\tilde{L})')'$ is available. As we will show,
the asymptotic properties
of the proposed one-step $R$-estimators will not be affected by the
choice of $\tilde\vartheta$. Practical choices will be provided in
Section \ref{simusec}.

Describing our one-step $R$-estimators requires:
\renewcommand{\theassumption}{(\Alph{assumption})}
\begin{assumption}\label{AssA}
For all $r\neq s\in\{1,\ldots,p\}$, we dispose of
sequences of estimators $\hat\gamma_{rs}(f)$ and $\hat\rho_{rs}(f)$
that: (i) are locally asymptotically discrete and that (ii), for
any $g\in{\mathcal{F}}_{\mathrm{ulan}}$, satisfy
$\hat\gamma_{rs}(f)=\gamma_{rs}(f,g)+o_\mathrm{P}(1)$
and
$\hat\rho_{rs}(f)=\rho_{rs}(f,g)+o_\mathrm{P}(1)$ as $\ny$,
under $\bigcup
_{\vartheta\in\Theta}\{\mathrm{P}^{(n)}_{\vartheta,g}\}$.
\end{assumption}

Sequences of estimators fulfilling this assumption will be provided in
Section~\ref{estimcrosssec} below. At this point, just note that
plugging\vspace*{2pt} in (\ref{definGstar}) the estimators from Assumption
\ref{AssA} and the\vspace*{-1pt} preliminary estimator $\tilde L$, defines a
statistic---$\hat\Gamma
^*_{\tilde{L},f;2}$, say---that 
consistently estimates $\Gamma^*_{{L},f,g;2}$ under $\bigcup_{\vartheta
\in
\Theta} \{\mathrm{P}^{(n)}_{\vartheta,g}\}$.

For any target density $f$, we propose the one-step $R$-estimator $\hat
{\underline{L}}_{f}$, with values in $\mathcal{M}_{1p}$, defined by
%
%
\begin{equation}\label{definpseudo}
\operatorname{vecd}^\circ\hat{\underline{L}}_{f}
:=
(\operatorname{vecd}^\circ\tilde{L})
+
n^{-1/2}
(\hat\Gamma^*_{\tilde{L},f;2})^{-1}
\underline{\Delta}^*_{\tilde{\vartheta},f;2}.
\end{equation}
%
The following result states the asymptotic properties of this estimator
(see Appendix \ref{appensecB} for a proof).
%
\begin{Theor}\label{theoestim}
Let Assumption \ref{AssA} hold, and fix
$f\in\underline{\mathcal{F}}_{\mathrm{ulan}}$.
Then
\textup{(i)}
under $\mathrm{P}^{(n)}_{\vartheta,g}$, with $\vartheta=(\mu',(
\operatorname{vecd}^\circ L)')'
\in\Theta$ and $g\in{\mathcal{F}}_{\mathrm{ulan}}$, we have that
%
%
%
%
\begin{eqnarray}\quad
\label{restim}
\sqrt{n} \operatorname{vec}(\hat{\underline{L}}_{f}-L)
&=&
C'(\Gamma^*_{L,f,g;2})^{-1}
\underline{\Delta}^*_{\vartheta,f;2}+o_\mathrm{P}(1)
\\
\label{represestim}
&=&
C'(\Gamma^*_{L,f,g;2})^{-1}
{\Delta}^*_{\vartheta,f,g;2}+o_\mathrm{P}(1)
\\
\label{estimasymplawvec} &\stackrel{\mathcal{L}}{\to}&
\mathcal{N}_{p(p-1)} ( 0, C'( \Gamma^*_{L,f,g;2})^{-1} \Gamma^*_{L,f;2}
(\Gamma^*_{L,f,g;2})^{-1\prime}C )
\end{eqnarray}
as $\ny$, where $\Delta^{*}_{\vartheta,f,g;2}$ is defined in
Theorem \ref{th2} (see Appendix \ref{appensec}).
\textup{(ii)}
The estimator $\hat{\underline{L}}_{f}$ is semiparametrically
efficient at $f$.
\end{Theor}

The result in (\ref{restim}) justifies calling $\hat{\underline
{L}}_{f}$ an $R$-estimator
since it shows that $n^{1/2}(\hat{\underline{L}}_{f}-L)$ is
asymptotically equivalent to a random matrix that is measurable with
respect to the signed ranks $S_i(\vartheta),R^{+}_i(\vartheta)$
in (\ref
{signedranks}). The asymptotic equivalence in (\ref{represestim}) gives
a Bahadur-type representation result for $\hat{\underline{L}}_{f}$ with
summands that are independent and identically distributed, hence leads
trivially to the
asymptotic normality result in (\ref{estimasymplawvec}). Recalling
that\vspace*{-1pt} $\hat\Gamma^*_{\tilde{L},f;2}$ consistently estimates $\Gamma
^*_{L,f,g;2}$ under $\bigcup_{\vartheta\in\Theta} \{\mathrm
{P}^{(n)}_{\vartheta
,g}\}$, it is\vspace*{1pt} clear that asymptotic (signed-rank) confidence zones
for $L$ may easily be obtained from this asymptotic normality result.

For $r\neq s\in\{1,\ldots,p\}$, define $\hat\alpha_{rs}(f)$
and $\hat
\beta_{rs}(f)$ as the statistics obtained by plugging the
estimators $\hat\gamma_{rs}(f)$ and $\hat\rho_{rs}(f)$ from
Assumption \ref{AssA} in
%
%
\begin{equation}
\label{alphabeta}
\cases{
\displaystyle \alpha_{rs}(f,g)
:=
\frac{\gamma_{rs}(f,g)}{\gamma_{rs}(f,g)\gamma_{sr}(f,g)-\rho
_{rs}(f,g)\rho_{sr}(f,g)}
\vspace*{6pt}\cr
\displaystyle \beta_{rs}(f,g)
:=
\frac{-\rho_{rs}(f,g)}{\gamma_{rs}(f,g)\gamma_{sr}(f,g)-\rho
_{rs}(f,g)\rho_{sr}(f,g)},
}
\end{equation}
and let $\hat\alpha_{rr}(f):=0=:\hat\beta_{rr}(f)$, $r=1,\ldots,p$.
The estimator $\hat{\underline{L}}_{f}$ then admits the following
explicit expression (see Appendix \ref{appensecB} for a proof).
\begin{Theor}\label{theoestimexplicit}
Let Assumption \ref{AssA} hold, and fix $f\in\underline{\mathcal{F}}_{\mathrm{ulan}}$.
Let $\hat N_{f}:=(\hat{\mathcal{A}}_{f}' \odot\underline
{T}_{\tilde
\vartheta,f}) + (\hat{\mathcal{B}}_{f}'\odot\underline
{T}_{\tilde
\vartheta,f}')$, where we let $\hat{\mathcal{A}}_{f}:=(\hat\alpha
_{rs}(f))$ and $\hat{\mathcal{B}}_{f}:=(\hat\beta_{rs}(f))$.
Then
the estimator $\hat{\underline{L}}_{f}$ rewrites
%
%
\begin{equation}\label{finalestim}
\hat{\underline{L}}_{f} = \tilde{L} + \frac{1}{\sqrt{n}} \tilde{L} [
\hat N_{f} - \operatorname{diag}(\tilde{L} \hat N_{f}) ] ,
\end{equation}
where $\operatorname{diag}(A)=A-\operatorname{odiag}(A)$ stands for the
diagonal matrix
with the same diagonal entries as $A$.
\end{Theor}

It is
straightforward\vspace*{-1pt} to check that the role of the term $-\frac{1}{\sqrt
{n}}\tilde{L} \operatorname{diag}(\tilde{L}\hat N_{f})$ in the
one-step correction $\frac{1}{\sqrt{n}}%
\tilde{L} [ \hat N_{f} -\operatorname{diag}(\tilde{L}\hat N_{f}) ]$ of
$\tilde{L}$
is merely\vspace*{-1pt} to ensure that the diagonal entries of $\hat{\underline
{L}}_{f}$ remain equal to one, hence that $\hat{\underline{L}}_{f}$
takes values in $\mathcal{M}_{1p}$ (for $n$ large enough).


As shown above, the estimator $\hat{\underline{L}}_{f}$ enjoys very
nice properties: its asymptotic behavior is completely characterized,
it is semiparametrically efficient under correctly specified densities,
yet remains root-$n$ consistent and asymptotically normal under a broad
range of densities $g$, its asymptotic covariance matrix can easily be
estimated consistently, etc.

However,\vspace*{1pt} $\hat{\underline{L}}_{f}$ requires estimates $\hat\gamma
_{rs}(f)$ and $\hat\rho_{rs}(f)$ that fulfill Assumption \ref{AssA}. We now
provide such estimates.

\subsection{Estimation of cross-information coefficients}
\label{estimcrosssec}

Of course, it is always possible to estimate consistently the
cross-informa\-tion coefficients $\gamma_{rs}(f,g)$ and $\rho
_{rs}(f,g)$ by replacing $g$ in (\ref{crossinfocoeffgam}) and (\ref
{crossinfocoeffrho}) with appropriate window or kernel density
estimates---this can be achieved since the residuals $Z_{ir}(\tilde
\vartheta)$, $i=1,\ldots,n$ typically are asymptotically i.i.d.
with density $g_r$. Rank-based methods, however, intend to
eliminate---through invariance arguments---the nuisance $g$ without
estimating it, so that density estimation methods simply are antinomic
to the spirit of rank-based methods. 

Therefore, we rather propose a solution that is based on ranks and
avoids estimating the underlying nuisance $g$. The method, that relies
on the asymptotic linearity---under $g$---of an appropriate rank-based
statistic $\underline{S}_{\vartheta,f}$, was first used in~\cite
{Ha06c}, where there is only one cross-information coefficient $J(f,g)$
to be estimated. There, it is crucial that $J(f,g)$ is involved as a
scalar factor in the asymptotic covariance matrix, under $g$, between
the rank-based efficient central sequence $\underline{\Delta
}^*_{\vartheta,f}$ and the parametric central sequence ${\Delta
}_{\vartheta,g}$. In \cite{r12}, the method was extended to allow for
the estimation of a cross-information coefficient that appears as a
scalar factor in the linear term of the asymptotic linearity,
under $g$, of a (possibly vector-valued) rank-based
statistic $\underline{S}_{\vartheta,f}$.

In all cases, thus, this method was only used to estimate a \textit
{single} cross-infor\-mation coefficient that appears as a \textit{scalar
factor} in some structural---typically, cross-information---matrix. In
this respect, our problem, which requires us to estimate $2p(p-1)$
cross-information quantities appearing in various entries of the
cross-information matrix $\Gamma^*_{L,f,g;2}$, is much more complex.
Yet, as we now show, it allows for a solution relying on the same basic
idea of exploiting the asymptotic linearity, under $g$, of an
appropriate $f$-score rank-based statistic.

Based on the preliminary estimator $ \tilde\vartheta := ( \tilde\mu',
(\operatorname{vecd}^\circ\tilde{L})')' $ at hand, define
$ \tilde\vartheta_{\lambda}^{\gamma_{rs}} := ( \tilde\mu',
(\operatorname{vecd}^\circ\tilde{L}_{\lambda}^{\gamma_{rs}})')'$,
$\lambda\geq0$, with
\[
%
%
\tilde{L}_{\lambda}^{\gamma_{rs}}:= \tilde{L} +n^{-1/2}\lambda
(\underline{T}_{\tilde\vartheta,f})_{rs} \tilde{L} \bigl( e_re_s'
-\operatorname{diag}(\tilde{L} e_re_s') \bigr),
\]
and $ \tilde\vartheta_{\lambda}^{\rho_{rs}} := ( \tilde\mu',
(\operatorname{vecd}^\circ\tilde{L}_{\lambda}^{\rho_{rs}})')' $,
$\lambda\geq0$, with
\[
%
%
\tilde{L}_{\lambda}^{\rho_{rs}}:= \tilde{L} +n^{-1/2}\lambda
(\underline{T}_{\tilde\vartheta,f})_{sr} \tilde{L} \bigl( e_re_s' -
\operatorname{diag}(\tilde{L} e_re_s') \bigr);
\]
note that, at $\lambda=0$, $\tilde\vartheta_{\lambda}^{\gamma
_{rs}}=\tilde\vartheta_{\lambda}^{\rho_{rs}}=\tilde\vartheta$.
We then have the following result that is crucial for the construction
of the estimators $\hat\gamma_{rs}(f)$ and $\hat\rho_{rs}(f)$; see
Appendix \ref{appensecB} for a proof.
\begin{Lem}\label{gamrhoestimtheor}
Fix $\vartheta\in\Theta$, $f\in\underline{\mathcal{F}}_{\mathrm{ulan}}$,
$g\in{\mathcal{F}}_{\mathrm{ulan}}$ and $r\neq s\in\{1,\ldots,p\}$. Then
$
h^{\gamma_{rs}}(\lambda)
:=
(\underline{T}_{\tilde\vartheta,f})_{rs}
(\underline{T}_{\tilde\vartheta_{\lambda}^{\gamma_{rs}},f})_{rs}
=
(1-\lambda\gamma_{rs}(f,g))
((\underline{T}_{\tilde\vartheta,f})_{rs})^2
+o_\mathrm{P}(1)
$
and
$
h^{\rho_{rs}}(\lambda)
:=
(\underline{T}_{\tilde\vartheta,f})_{sr}
(\underline{T}_{\tilde\vartheta_{\lambda}^{\rho_{rs}},f})_{sr}
=
(1-\lambda\rho_{rs}(f,g))
((\underline{T}_{\tilde\vartheta,f})_{sr})^2
+o_\mathrm{P}(1)
$
as $\ny$, under $\mathrm{P}^{(n)}_{\vartheta,g}$.
\end{Lem}

The mappings $\lambda\mapsto h^{\gamma_{rs}}(\lambda)$ and $\lambda
\mapsto h^{\rho_{rs}}(\lambda)$ assume a positive value in $\lambda=0$,
and, as shown by Lemma \ref{gamrhoestimtheor}, are---up to $o_\mathrm
{P}(1)$'s as $\ny$ under $\mathrm{P}^{(n)}_{\vartheta,g}$---monotone decreasing
functions that become negative at $\lambda=(\gamma_{rs}(f,g))^{-1}$ and
$\lambda=(\rho_{rs}(f,g))^{-1}$, respectively. Restricting to a grid of
values of the form $\lambda_j=j/c$ for some large discretization
constant $c$ (which is needed to achieve the required discreteness),
this naturally leads---via linear interpolation---to the
estimators $\hat\gamma_{rs}(f)$ and $\hat\rho_{rs}(f)$ defined through
%
%
\begin{eqnarray}
\label{gammaestimdiscr}
(\hat\gamma_{rs}(f))^{-1}
:\!&=&
\lambda_{\gamma_{rs}}:=
\lambda^{-}_{\gamma_{rs}}
+\frac
{(\lambda^{+}_{\gamma_{rs}}-\lambda^{-}_{\gamma_{rs}})h^{\gamma
_{rs}}(\lambda^{-}_{\gamma_{rs}})}
{h^{\gamma_{rs}}(\lambda^{-}_{\gamma_{rs}})-h^{\gamma_{rs}}(\lambda
^{+}_{\gamma_{rs}})}
\nonumber\\[-8pt]\\[-8pt]
&=& \lambda^{-}_{\gamma_{rs}} + \frac {c^{-1}
h^{\gamma_{rs}}(\lambda^{-}_{\gamma_{rs}})}
{h^{\gamma_{rs}}(\lambda^{-}_{\gamma_{rs}})-h^{\gamma_{rs}}(\lambda
^{+}_{\gamma_{rs}})} \nonumber
\end{eqnarray}
with
$
\lambda^{-}_{\gamma_{rs}}
:=
\inf
\{j\in\N\dvtx
h^{\gamma_{rs}}(\lambda_{j+1})<0
\}
$
and
$
\lambda^{+}_{\gamma_{rs}}:=
\lambda^{-}_{\gamma_{rs}}+\frac{1}{c}
$,
and
%
%
\begin{eqnarray}
\label{rhoestimdiscr}
(\hat\rho_{rs}(f))^{-1} := \lambda_{\rho_{rs}}
&:=&
\lambda^{-}_{\rho_{rs}} + \frac {c^{-1}
h^{\rho_{rs}}(\lambda^{-}_{\rho_{rs}})}
{h^{\rho_{rs}}(\lambda^{-}_{\rho_{rs}})-h^{\rho_{rs}}(\lambda
^{+}_{\rho_{rs}})}
\end{eqnarray}
with\vspace*{1pt} $ \lambda^{-}_{\rho_{rs}} := \inf \{j\in\N\dvtx
h^{\rho_{rs}}(\lambda_{j+1})<0 \} $ and $ \lambda^{+}_{\rho_{rs}}:=
\lambda^{-}_{\rho_{rs}}+\frac{1}{c} $.
We have the following result
(see the supplemental article \cite{IP11a} for a proof).
\begin{Theor} \label{gamrhoestimconsist}
Fix $\vartheta\in\Theta$,
$f\in\underline{\mathcal{F}}_{\mathrm{ulan}}$, and
$g\in{\mathcal{F}}_{\mathrm{ulan}}$. Assume
that $\tilde\vartheta$ is such that, for all $\varepsilon>0$, there
exist $\delta_{\varepsilon}>0$ and an integer $N_{\varepsilon}$ such that
%
%
\begin{equation}
\label{assB} \mathrm{P}^{(n)}_{\vartheta,g} [
(\underline{T}_{\tilde\vartheta,f})_{rs}\geq\delta_{\varepsilon}
]\geq1-\varepsilon
\end{equation}
for all $n\geq N_{\varepsilon}$, $r\neq s\in\{1,\ldots,p\}$.
Then, for any such $r,s$, $\hat\gamma_{rs}(f)
=\gamma_{rs}(f,g)+o_\mathrm{P}(1)
$
and
$\hat\rho_{rs}(f)
=\rho_{rs}(f,g)+o_\mathrm{P}(1)
$, as $\ny$ under $\mathrm{P}^{(n)}_{\vartheta,g}$, hence $\hat\gamma_{rs}(f)$
and
$\hat\rho_{rs}(f)$ satisfy Assumption \ref{AssA}.
\end{Theor}


We point out that the assumption in (\ref{assB}) is extremely mild, as
it only requires that there is
no couple $(r,s)$, $r\neq s$, for which $(\underline{T}_{\tilde
\vartheta
,f})_{rs}$ asymptotically has an atom in zero. It therefore rules out
preliminary estimators $\tilde L$ defined through the (rank-based)
$f$-likelihood equation $(\underline{T}_{\vartheta,f})_{rs}=0$. 



\section{Simulations}
\label{simusec}

Here we report simulation results for point estimation
only---simulation results for hypothesis testing can be found in the
supplemental article~\cite{IP11a}. Our aim is to both compare the
proposed estimators with some competitors and to investigate the
validity of asymptotic results.

We used the following competitors: (i) FastICA from \cite{Hy97,Hy99},
which is by far the most commonly used estimate in practice; we used
here its deflation based version with the standard nonlinearity
function pow3. (ii) FOBI from \cite{Ca89}, which is one of the earliest
solutions to the ICA problem and is often used as a benchmark estimate.
(iii) The estimate based on two scatter matrices from \cite{SCA}; here
the two scatter matrices used are the regular empirical covariance
matrix (COV) and the van der Waerden rank-based estimator (HOP)
from \cite{Ha06c} (actually, HOP is not a scatter matrix but rather a
shape matrix, which is allowed in \cite{SCA}). Root-$n$ consistency of
the resulting estimates $\hat L_{\mathrm{FICA}}$, $\hat L_{\mathrm{FOBI}}$
and $\hat L_{\mathrm{COV}\_\mathrm{HOP}}$ of $L$ requires finite
sixth-, eighth- and
fourth-order moments, respectively, and follows from \cite{Il10,Il11}
and \cite{Ol10}.

%

We focused on the bivariate case $p=2$, and we generated, for three
different setups indexed by $d\in\{1,2,3\}$, $M=2\mbox{,}000$ independent
random samples $Z^{(d,m)}_i=(Z^{(d,m)}_{i1},Z^{(d,m)}_{i2})'$,
$i=1,\ldots,n$, of size $n=4\mbox{,}000$. Denoting
by $g^{(d)}(z)=g^{(d)}_1(z_1)g^{(d)}_2(z_2)$ the common pdf of
$Z^{(d,m)}_i$, $i=1,\ldots,n$, $m=1,\ldots,M$, the marginal
densities $g^{(d)}_1$ and $g^{(d)}_2$ were chosen as follows:
\begin{longlist}
\item
In Setup $d=1$, $g^{(d)}_1$ is the pdf of the standard normal
distribution ($\mathcal{N}$), and $g^{(d)}_2$ is the pdf of the Student
distribution with $5$ degrees of freedom ($t_5$);
\item
In Setup $d=2$, $g^{(d)}_1$ is the pdf of the logistic distribution
with scale parameter one (log), and $g^{(d)}_2$ is $t_5$;
\item
%
In Setup $d=3$, $g^{(d)}_1$ is $t_8$ and $g^{(d)}_2$ is $t_5$.
\end{longlist}
We chose to use $L=I_2$ and $\mu=(0,0)'$, so that the observations are
given by $X_i^{(d,m)}= LZ^{(d,m)}_i+\mu=Z^{(d,m)}_i$ (other values
of $L$ and $\mu$ led to extremely similar results).

For each sample, we computed the competing estimates $\hat L_{
\mathrm{FICA}}$, $\hat L_{\mathrm{FOBI}}$ and $\hat
L_{\mathrm{COV}\_\mathrm{HOP}}$ defined above.
Each of these were also used as a preliminary estimator $\tilde L$ in the
construction of three $R$-estimators: $\hat{\underline{L}}_{f^{(j)}}$,
$j=1,2,3$,
with $f^{(j)}=g^{(j)}$ for all $j$. In the resulting nine $R$-estimators,
we used the location estimate $\hat\mu=\tilde L \operatorname
{Med}[\tilde
L^{-1}X_1,\ldots,\tilde L^{-1}X_n]$, based on the preliminary
estimate $\tilde L$ used to initiate the one-step procedure.

\begin{figure}

\includegraphics{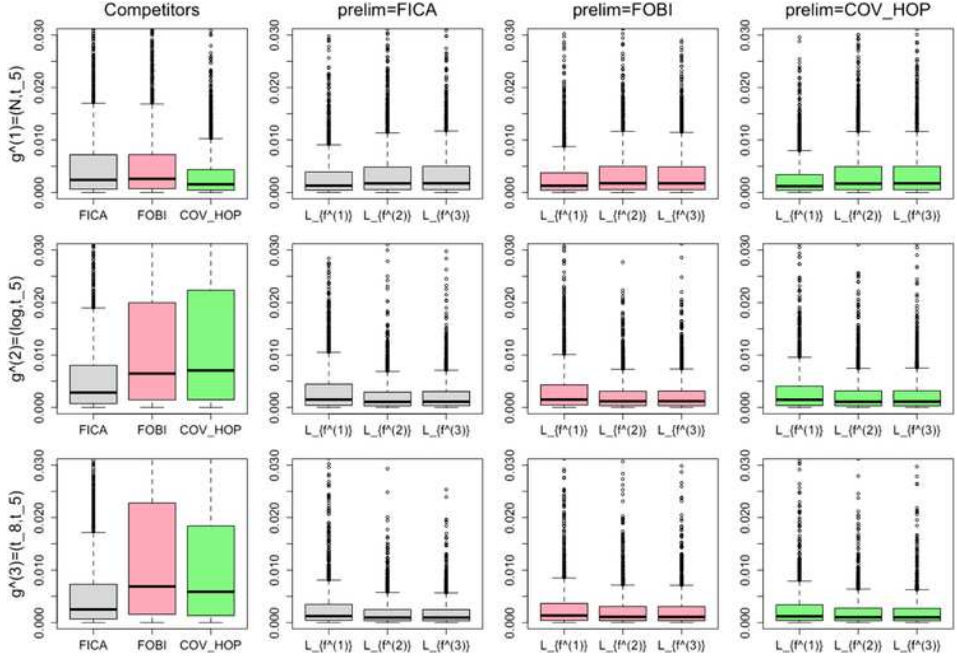}

\caption{Boxplots of the squared errors $\| \hat L-L \|^2$ [see
(\protect\ref{errorfr})] obtained in $M=2\mbox{,}000$ replications from
setups $d=1,2,3$ (associated\vspace*{1pt} with underlying distributions $g^{(d)}$,
$d=1,2,3$) for the competitors $\hat{L}_{\mathrm{FICA}}$, $\hat{L}_{
\mathrm{FOBI}}$ and $\hat{L}_{\mathrm{COV}\_\mathrm{HOP}}$, and the
nine $R$-estimators $\hat L_f$ resulting from all combinations of a
target density $f^{(j)}=g^{(j)}$, $j=1,2,3$, and one of the three
preliminary estimators $\hat{L}_{\mathrm{FICA}}$,
$\hat{L}_{\mathrm{FOBI}}$ and $\hat {L}_{\mathrm{COV}\_\mathrm{HOP}}$;
see Section \protect\ref{simusec} for details. The sample size is
$n=4\mbox{,}000$.} \label{fig2revn4000}
\end{figure}

Figure \ref{fig2revn4000} reports, for each setup $d$, a boxplot of
the $M$ squared errors
%
%
\begin{equation} \label{errorfr}\qquad
\bigl\| \hat L\bigl(X_1^{(d,m)},\ldots,X_n^{(d,m)}\bigr)-L \bigr\|^2
= \sum_{\stackrel{r,s=1}{r\neq s}}^p \bigl( \hat
L_{rs}\bigl(X_1^{(d,m)},\ldots,X_n^{(d,m)}\bigr)-L_{rs} \bigr)^2
\end{equation}
for each of the twelve estimators $\hat L$ considered (the nine
$R$-estimators and their three competitors).

%
%

The results show that, in each setup, all $R$-estimators dramatically
improve over their competitors.
The behavior of the $R$-estimators does not much depend on the
preliminary estimator $\tilde{L}$ used. Optimality of $\hat
{L}_{f^{(d)}}$ in Setup $d$ is confirmed. Most importantly, as stated
for hypothesis testing at the end of Section \ref{testsec}, the
performances of the $R$-estimators do not depend much on the target
density $f^{(j)}$ adopted, so that one should not worry much about the
choice of the target density in practice. Quite surprisingly,
$R$-estimators behave remarkably well even when based on preliminary
estimators that, due to heavy tails, fail to be root-$n$ consistent.

\begin{figure}

\includegraphics{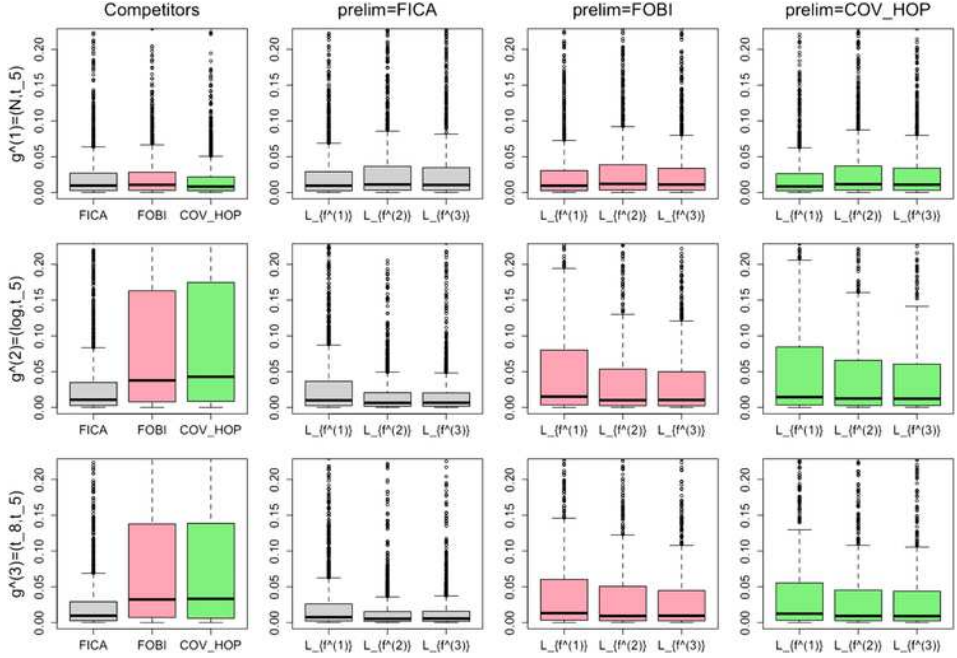}

\caption{The same boxplots as in Figure \protect\ref{fig2revn4000},
but based on sample size $n=800$.}
\label{fig2revn800}
\end{figure}

In order to investigate small-sample behavior of the estimates, we
reran the exact same simulation with sample size $n=800$; in ICA, where
most applications involve sample sizes that are not in hundreds, but
much larger, this sample size can indeed be considered small. Results
are reported in Figure \ref{fig2revn800}. They indicate that, in
Setups 2 and 3, $R$-estimators still improve significantly over their
competitors, and particularly over $\hat{L}_{\mathrm{FOBI}}$ and $\hat
{L}_{\mathrm{COV}\_\mathrm{HOP}}$.
In Setup 1, there seem to be no improvement. Compared to results
for $n=4\mbox{,}000$, the behavior of one-step $R$-estimators here depends more
on the preliminary estimator used. Performances of $R$-estimators again
do not depend crucially on the target density, and optimality under
correctly specified densities is preserved in most cases.

%
%

As a conclusion, for practical sample sizes, the proposed $R$-estimators
outperform the standard competitors considered, and their behavior is
very well in line with our asymptotic results.

\begin{figure}

\includegraphics{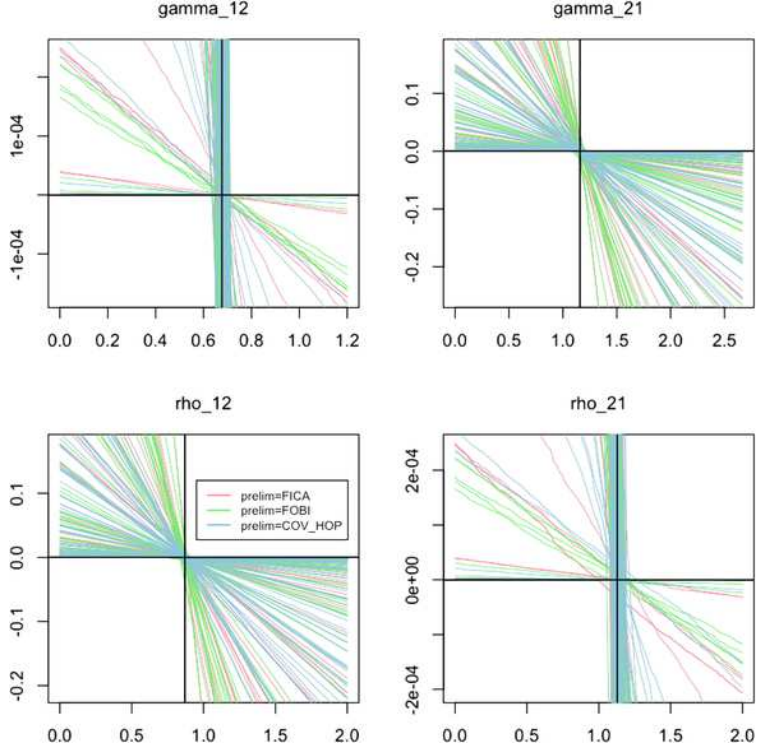}

\caption{Top left: 150 graphs of the mapping $\lambda\mapsto
h^{\gamma_{12}}(\lambda)$ based on $f=f^{(3)}$, associated with the
first 50 replications from Setup 1 ($g=g^{(1)} $) in Figure
\protect\ref{fig2revn4000} (sample size is $n=4\mbox{,}000$): the 50
curves in pink, green, and blue are based on the preliminary estimators
$\hat {L}_{\mathrm{FICA}}$, $\hat{L}_{\mathrm{FOBI}}$ and
$\hat{L}_{\mathrm {COV}\_\mathrm{HOP}}$, respectively. Top right,
bottom left, and bottom right: the corresponding plots for the mappings
$\lambda\mapsto h^{\gamma_{21}}(\lambda)$, $\lambda\mapsto
h^{\rho_{12}}(\lambda )$ and $\lambda\mapsto h^{\rho_{21}}(\lambda)$,
respectively.} \label{fig3rev}
\end{figure}

Finally, we illustrate the proposed method for estimating
cross-information coefficients. We consider again the first 50
replications of our simulation with $n=4\mbox{,}000$, and focus on
Setup 1 ($g=g^{(1)}$) and the target density $f=f^{(3)}$ (\mbox{$\neq$}$g^{(1)})$.
The cross-information coefficients to be estimated then are
$\gamma_{12}(f,g)\approx1.478$, %
$\gamma_{21}(f,g)\approx0.862$,
$\rho_{12}(f,g)\approx1.149$
and
$\rho_{21}(f,g)\approx0.887$. The upper left picture in Figure
\ref{fig3rev} shows 150 graphs of the mapping $\lambda\mapsto
h^{\gamma_{12}}(\lambda)$ (based on $f=f^{(3)}$), among which the 50
%
%
pink curves are based on $\tilde{L}=\hat{L}_{\mathrm{FICA}}$, the 50
green\vspace*{1pt} curves are based on
$\tilde{L}=\hat{L}_{\mathrm{FOBI}}$, and the 50 blue ones are based on
$\tilde{L}=\hat{L}_{\mathrm{COV}\_ \mathrm{HOP}}$.
The upper right, bottom left and bottom right pictures of the same
figure provide the corresponding graphs for the mappings
$\lambda\mapsto h^{\gamma_{21}}(\lambda)$,
$\lambda\mapsto h^{\rho_{12}}(\lambda)$,
and
$\lambda\mapsto h^{\rho_{21}}(\lambda)$, respectively. The value at
which each graph crosses the $\lambda$-axis is the resulting estimate
of the inverse of the associated cross-information coefficient. To be
able to evaluate the results, we plotted, in each picture, a vertical
black line at the corresponding theoretical value, namely at $1/\gamma
_{12}(f,g)$, %
$1/\gamma_{21}(f,g)$,
$1/\rho_{12}(f,g)$
and
$1/\rho_{21}(f,g)$. Clearly, the results are excellent, and there does
not seem to be much dependence on the preliminary estimator $\tilde L$ used.


\begin{appendix}\label{app}

\section{Rank-based efficient central sequences}
\label{appensec}

In this first Appendix, we study the asymptotic behavior of the
rank-based efficient central sequences $\underline{\Delta
}^{*}_{\vartheta,f;2}$. The main result is the following (see
Appendix \ref{appensecB} for a proof).
%
%
\begin{Theor}\label{th2}
Fix $\vartheta=(\mu',(\operatorname{vecd}^\circ L)')'\in\Theta$ and
$f\in
\underline{\mathcal{F}}_{\mathrm{ulan}}$. Then, \textup{(i)} for any $g\in\mathcal{F}$,
\[
\underline{\Delta}^{*}_{\vartheta,f;2}
=
\Delta^{*}_{\vartheta,f,g;2}
+o_{L^2}(1)
\]
as $\ny$, under $\mathrm{P}^{(n)}_{\vartheta,g}$,
where
$
\Delta^{*}_{\vartheta,f,g;2}
:=
C (I_p\otimes L^{-1})'
\operatorname{vec}[\operatorname{odiag}(\frac{1}{\sqrt{n}}\sum_{i=1}^n
(S_i\odot\varphi_{f}(F_+^{-1}(G_+(|Z_{i}|))))(S_i\odot
F_+^{-1}(G_+(|Z_{i}|)))')].
$
\textup{(ii)} Under $\mathrm{P}^{(n)}_{\vartheta+n^{-1/2}\tau,g}$, with $\tau
=(\tau
_1',\tau_2')' \in\R^p\times\R^{p(p-1)}$ and $g\in{\mathcal
{F}}_{\mathrm{ulan}}$,
\[
\underline{\Delta}^{*}_{\vartheta,f;2}
\stackrel{\mathcal{L}}{\to}
\mathcal{N}_{p(p-1)}( \Gamma^*_{L,f,g;2}\tau_2, \Gamma^*_{L,f;2})
\]
as $\ny$ (for $\tau=0$, the result only requires that $g\in\mathcal{F}$).
\textup{(iii)}
Still with $\tau=(\tau_1',\tau_2')' \in\R^p\times\R^{p(p-1)}$ and
$g\in{\mathcal{F}}_{\mathrm{ulan}}$,
$
\underline{\Delta}^{*}_{\vartheta+n^{-1/2}\tau,f;2}-\underline
{\Delta
}_{\vartheta,f;2}^{*}
=
-\Gamma^*_{L,f,g;2} \tau_2
+o_\mathrm{P}(1)
$
as $\ny$, under $\mathrm{P}^{(n)}_{\vartheta,g}$. 
\end{Theor}

Both
for hypothesis testing and point estimation, we had to replace in
$\underline{\Delta}^{*}_{\vartheta,f;2}$ the parameter $\vartheta$ with
some estimator ($\check\vartheta^{(n)}$, say). The asymptotic
behavior of
the resulting (so-called \textit{aligned}) rank-based efficient central
sequence $\underline{\Delta}^{*}_{\check\vartheta^{(n)},f;2}$ is
given in
the following result.
\begin{Corol}\label{corol}
Fix $\vartheta=(\mu',(\operatorname{vecd}^\circ L)')'\in\Theta$,
$f\in
\underline{\mathcal{F}}_{\mathrm{ulan}}$, and $g\in{\mathcal
{F}}_{\mathrm{ulan}}$. Let $\check\vartheta=\check\vartheta
^{(n)}=(\check\mu
',(\operatorname{vecd}^\circ\check L)')'$ be a locally asymptotically discrete
sequence of random
vectors satisfying $n^{1/2}(\check\vartheta-\vartheta)=O_\mathrm{P}(1)$
as $\ny$, under $\mathrm{P}^{(n)}_{\vartheta,g}$. Then
$
\underline{\Delta}_{\check\vartheta,f;2}^{*}-\underline{\Delta
}_{\vartheta,f;2}^{*}
=
-\Gamma^*_{L,f,g;2}n^{1/2}\operatorname{vecd}^\circ(\check L-L)
+o_\mathrm{P}(1),
$
still as $\ny$, under~$\mathrm{P}^{(n)}_{\vartheta,g}$.
\end{Corol}
%


Since the\vspace*{1pt} sequence of estimators $\check\vartheta^{(n)}$ is assumed to
be locally asymptotically discrete [which means that the number of
possible values of $\check\vartheta^{(n)}$ in balls with $O(n^{-1/2})$
radius centered at $\vartheta$ is bounded as $\ny$], this result is a
direct consequence of Theorem \ref{th2}(iii) and Lemma 4.4 from
\cite{kreiss}.
Local asymptotic discreteness is a concept that goes back to Le Cam and
is quite standard in one-step estimation; see,\vspace*{1pt} for example, \cite{Bi82}
or \cite{kreiss}.

Of course,\vspace*{1pt} a sequence of estimators
$\check\vartheta^{(n)}$ can always be discretized by replacing each
component $(\check\vartheta^{(n)})_\ell $ with
\[
\bigl(\check\vartheta^{(n)}_\#\bigr)_\ell:=(cn^{1/2})^{-1}\operatorname
{sign}\bigl(\bigl(\check \vartheta ^{(n)}\bigr)_\ell\bigr)\bigl\lceil
cn^{1/2}\bigl\vert\bigl(\check\vartheta^{(n)}\bigr)_\ell\bigr\vert
\bigr\rceil,\qquad \ell=1,\ldots,p^2,
\]
for some arbitrary constant $c>0$.
In practice, however, one can safely forget about such discretizations:
irrespective of the accuracy of the computer used, the discretization
constant $c$ can always be chosen large enough to make discretization
be irrelevant at the fixed sample size $n_0$ at hand---hence also at
any $n> n_0$.

\section{Proofs}
\label{appensecB}

\subsection{\texorpdfstring{Proofs of Theorems \protect\ref{th1} and
\protect\ref{th2}}{Proofs of Theorems 2.1 and A.1}}

The proofs of this section make use of the H\'{a}jek projection theorem
for linear signed-rank statistics (see, e.g., \cite{Pu85}, Chapter 3),
which states that, if $Y_i=\operatorname{Sign}(Y_i)|Y_i|$, $i=1,\ldots
,n$, are
i.i.d. with (absolutely continuous) cdf $G$ and if $K\dvtx(0,1)\to
\R
$ is a continuous and square-integrable score function that can be
written as the difference of two monotone increasing functions, then
%
%
\begin{eqnarray}
&&\frac{1}{\sqrt{n}} \sum_{i=1}^n \operatorname{Sign}(Y_i) K(G_+(|Y_i|))
\nonumber
\\
\label{apprsco}
&&\qquad
=\frac{1}{\sqrt{n}} \sum_{i=1}^n \operatorname{Sign}(Y_i) K\biggl(\frac
{R^+_{i}}{n+1}\biggr)+o_{L^2}(1)
\\
\label{exasco2222}
&&\qquad
=\frac{1}{\sqrt{n}} \sum_{i=1}^n \operatorname{Sign}(Y_i) \mathrm{E}[
K(G_+(|Y_i|)) | R^+_{i}
]+o_{L^2}(1)
\end{eqnarray}
as $\ny$, where $G_+$ stands for the common cdf of the $|Y_i|$'s and
$R^+_{i}$ denotes the rank of $|Y_i|$ among $|Y_1|,\ldots,|Y_n|$.
The quantities in (\ref{apprsco})
and (\ref{exasco2222}) are linear signed-rank quantities that are said
to be based on \textit{approximate} and \textit{exact} scores, respectively.

In the rest of this section, we fix $\vartheta\in\Theta$, $f\in
\underline{\mathcal
{F}}_{\mathrm{ulan}}$, and $g\in\mathcal{F}$. We write throughout $Z_{i}$,
$S_{i}$, and $R^+_{i}$, for $Z_{i}(\vartheta)$, $S_{i}(\vartheta)$, and
$R^+_{i}(\vartheta)$, respectively. We also write $\mathrm{E}_{h}$ instead
of $\mathrm{E}^{(n)}_{\vartheta,h}$, with $h=f,g$.
We then start with the proof of Theorem~\ref{th2}(i). 
\begin{pf*}{Proof of Theorem \ref{th2}\textup{(i)}}
Fix $r\neq s\in\{1,\ldots,p\}$ and two score
functions $K_a,K_b\dvtx(0,1)\to\R$ with the same properties as $K$ above.
Then, by using (i) $\mathrm{E}_g[S_{ir}]=0$, (ii) the independence (under
$\mathrm{P}^{(n)}_{\vartheta,g}$) between the\vspace*{1pt} $S_{ir}$'s and the
$(R_{ir},|Z_{ir}|)$'s, and (iii) the independence between the
$Z_{ir}$'s and the $Z_{is}$'s, we obtain
%
\begin{eqnarray*}
&& \mathrm{E}_{g} \Biggl[ \Biggl( \frac{1}{\sqrt{n}} \sum_{i=1}^n
S_{ir}S_{is} \biggl( K_a(G_{+r}(|Z_{ir}|)) K_b(G_{+s}(|Z_{is}|))\\
&&\qquad\quad\hspace*{82.2pt}{}-K_a\biggl(\frac{R^+_{ir}}{n+1}\biggr) K_b\biggl(\frac{R^+_{is}}{n+1}\biggr)
\biggr) \Biggr)^2 \Biggr]
\\
&&\qquad = \frac{1}{n} \sum_{i=1}^n \mathrm{E}_{g} \biggl[ \biggl(
K_a(G_{+r}(|Z_{ir}|)) K_b(G_{+s}(|Z_{is}|))\\
&&\qquad\quad\hspace*{72.2pt}{}-K_a\biggl(\frac{R^+_{ir}}{n+1}\biggr) K_b\biggl(\frac{R^+_{is}}{n+1}\biggr)
\biggr)^2 \biggr]
\\
&&\qquad \leq2\mathrm{E}_{g} \biggl[ \biggl( K_a(G_{+r}(|Z_{ir}|))
-K_a\biggl(\frac{R^+_{ir}}{n+1}\biggr) \biggr)^2 \biggr] \mathrm{E}_{g}
[ K^2_b(G_{+s}(|Z_{is}|)) ]
\\
&&\qquad\quad{} + 2\mathrm{E}_{g} \biggl[
K^2_a\biggl(\frac{R^+_{ir}}{n+1}\biggr)\biggr] \mathrm{E}_{g} \biggl[
\biggl( K_b(G_{+s}(|Z_{is}|)) -K_b\biggl(\frac{R^+_{is}}{n+1}\biggr)
\biggr)^2 \biggr].
\end{eqnarray*}
Consequently, the square integrability of $K_a$, $K_b$, and the
convergence to zero of both $\mathrm{E}_{g} [( K_a(G_{+r}(|Z_{ir}|))
-K_a(\frac{R^+_{ir}}{n+1}) )^2]$ and\vspace*{1pt} $\mathrm{E}_{g} [(
K_b(G_{+r}(|Z_{is}|)) -K_b(\frac{R^+_{is}}{n+1}) )^2]$ [which directly
follows from (\ref{apprsco})] entail
\begin{eqnarray*}
&&\frac{1}{\sqrt{n}} \sum_{i=1}^n S_{ir}S_{is} K_a(G_{+r}(|Z_{ir}|))
K_b(G_{+s}(|Z_{is}|))
\\
&&\qquad = \frac{1}{\sqrt{n}} \sum_{i=1}^n S_{ir} S_{is}
K_a\biggl(\frac{R^+_{ir}}{n+1}\biggr)
K_b\biggl(\frac{R^+_{is}}{n+1}\biggr)+o_{L^2}(1)
\end{eqnarray*}
as $\ny$, under $\mathrm{P}^{(n)}_{\vartheta,g}$. Theorem \ref{th2}(i) follows
by taking $K_a=\varphi_{f_r}\circ F_{+r}^{-1}$ and $K_b= F_{+s}^{-1}$.
\end{pf*}

We go on with the proof of Theorem \ref{th1}, for which it is
important to note that, by proceeding as in the proof of Theorem \ref
{th2}(i) but with (\ref{exasco2222}) instead of (\ref{apprsco}), we
further obtain that
%
%
\begin{eqnarray}\label{exadiff}
&&\frac{1}{\sqrt{n}} \sum_{i=1}^n S_{ir}S_{is}
K_a(G_{+r}(|Z_{ir}|))
K_b(G_{+s}(|Z_{is}|))
\nonumber\\
&&\qquad
=
\frac{1}{\sqrt{n}} \sum_{i=1}^n S_{ir} S_{is}
K_a\biggl(\frac{R^+_{ir}}{n+1}\biggr)
K_b\biggl(\frac{R^+_{is}}{n+1}\biggr)+o_{L^2}(1)
\nonumber\\[-8pt]\\[-8pt]
&&\qquad = \frac{1}{\sqrt{n}} \sum_{i=1}^n S_{ir} S_{is} \mathrm{E}_g[
K_a(G_{+r}(|Z_{ir}|)) | R_{+{ir}}]
\nonumber\\
&&\qquad\quad
\hphantom{\frac{1}{\sqrt{n}} \sum_{i=1}^n}
{}\times
\mathrm{E}_g[
K_b(G_{+s}(|Z_{is}|)) | R_{+{is}}
]+o_{L^2}(1),
\nonumber
\end{eqnarray}
still as $\ny$ under $\mathrm{P}^{(n)}_{\vartheta,g}$.
\begin{pf*}{Proof of Theorem \ref{th1}}
It is sufficient to prove Theorem \ref{th1}(i) only, since, as already
mentioned at the end of Section \ref{invarianceranksec}, Theorem \ref
{th1}(ii) follows from (\ref{HHWW03}) and Theorem \ref{th1}(i).
That is, we have to show that, for any $r,s\in\{1,\ldots,p\}$,
%
%
\begin{eqnarray} \label{toprovtheo1}
&&\mathrm{E}_{f} \Biggl[ \frac{1}{\sqrt{n}}
\sum_{i=1}^n\bigl(\varphi_{f}(Z_i)Z_i'-I_p\bigr)_{rs} | S_1, \ldots,
S_n, R^+_1, \ldots, R^+_n\Biggr]\nonumber\\[-8pt]\\[-8pt]
&&\qquad=(\underline{T}_{\vartheta,f})_{rs}+o_{L^2}(1)\nonumber
\end{eqnarray}
as $\ny$, under $\mathrm{P}^{(n)}_{\vartheta,f}$.
Now, the left-hand side of (\ref{toprovtheo1}) rewrites
%
%
\begin{eqnarray}\label{totos}\quad
&&\mathrm{E}_{f} \Biggl[ \frac{1}{\sqrt{n}}
\sum_{i=1}^n\bigl(\varphi_{f}(Z_i)Z_i'-I_p\bigr)_{rs} | S_1, \ldots,
S_n, R^+_1, \ldots, R^+_n\Biggr]
\nonumber
\\
&&\qquad = \frac{1}{\sqrt{n}} \sum_{i=1}^n \mathrm{E}_{f} [
S_{ir}S_{is} \varphi_{f}(|Z_{ir}|)|Z_{is}|-\delta_{rs} | S_1, \ldots,
S_n, R^+_1, \ldots, R^+_n]
\\
&&\qquad = \frac{1}{\sqrt{n}} \sum_{i=1}^n \bigl( S_{ir}S_{is}
\mathrm{E}_{f} [ \varphi_{f}(|Z_{ir}|)|Z_{is}| | R^+_{1r}, \ldots,
R^+_{nr}, R^+_{1s}, \ldots, R^+_{ns}
] -\delta_{rs} \bigr) .\nonumber
\end{eqnarray}
For $r\neq s$, this yields
\begin{eqnarray*}
&&\mathrm{E}_{f} \Biggl[ \frac{1}{\sqrt{n}}
\sum_{i=1}^n\bigl(\varphi_{f}(Z_i)Z_i'-I_p\bigr)_{rs} | S_1, \ldots, S_n, R^+_1,
\ldots, R^+_n\Biggr]
\\[-3pt]
&&\qquad = \frac{1}{\sqrt{n}} \sum_{i=1}^n S_{ir}S_{is} \mathrm{E}_{f}
[ \varphi_{f}(|Z_{ir}|) | R^+_{1r}, \ldots, R^+_{nr}] \mathrm{E}_{f} [
|Z_{is}| | R^+_{1s}, \ldots, R^+_{ns}]
\\[-3pt]
&&\qquad = \frac{1}{\sqrt{n}} \sum_{i=1}^n S_{ir} S_{is}
\varphi_{f_r}\biggl(F_{+r}^{-1}\biggl(\frac{R^+_{ir}}{n+1}\biggr)\biggr)
F_{+r}^{-1}\biggl(\frac{R^+_{is}}{n+1}\biggr)+o_{L^2}(1)
\\[-3pt]
&&\qquad = (\underline{T}_{\vartheta,f})_{rs}+o_{L^2}(1)
\end{eqnarray*}
as $\ny$, under $\mathrm{P}^{(n)}_{\vartheta,f}$, where we have
used (\ref
{exadiff}), still with $K_a=\varphi_{f_r}\circ F_{+r}^{-1}$ and $K_b=
F_{+s}^{-1}$, but this time at $g=f$. This establishes (\ref
{toprovtheo1}) for $r\neq s$. As for $r=s$, (\ref{totos}) now entails
[writing $K_{ab}(u):=\varphi_{f}(F_{+r}^{-1}(u))\times F_{+r}^{-1}(u)$
for all~$u$]
%
%
\begin{eqnarray}
&&
\mathrm{E}_{f} \Biggl[ \frac{1}{\sqrt{n}}
\sum_{i=1}^n\bigl(\varphi_{f}(Z_i)Z_i'-I_p\bigr)_{rs} | S_1, \ldots, S_n, R^+_1,
\ldots, R^+_n \Biggr] \nonumber
\\[-3pt]
&&\qquad = \Biggl( \frac{1}{\sqrt{n}} \sum_{i=1}^n \mathrm{E}_{f} [
\varphi_{f}(|Z_{ir}|)|Z_{ir}| | R^+_{1r}, \ldots, R^+_{nr} ] \Biggr)
-\sqrt{n} \nonumber
\\[-3pt]
&&\qquad = \mathrm{E}_{f} \Biggl[ \frac{1}{\sqrt{n}} \sum_{i=1}^n
K_{ab}(F_{+r}(|Z_{ir}|)) | R^+_{1r}, \ldots, R^+_{nr} \Biggr] -\sqrt{n}
\nonumber
\\[-3pt]
\label{repr} &&\qquad = \frac{1}{\sqrt{n}} \sum_{i=1}^n
K_{ab}\biggl(\frac{R^+_{i}}{n+1}\biggr) -\sqrt{n}+o_{L^2}(1)
\\[-3pt]
&&\qquad = \frac{1}{\sqrt{n}} \sum_{i=1}^n K_{ab}\biggl(\frac{i}{n+1}\biggr)
-\sqrt{n}+o_{L^2}(1) \nonumber
\\[-3pt]
\label{sq} &&\qquad = \sqrt{n} \int_0^1 K_{ab}(u) \,du - \sqrt{n}
+o_{L^2}(1)
\\[-3pt]
\label{pa} &&\qquad = o_{L^2}(1),
\end{eqnarray}
still as $\ny$, under $\mathrm{P}^{(n)}_{\vartheta,f}$,\vspace*{1pt}
where (\ref {repr}), (\ref{sq}) and (\ref{pa}) follow from the
H\'{a}jek projection theorem for linear \textit{rank} (not
\textit{signed-rank}) statistics (see, e.g., \cite{Pu85}, Chapter 2),
the square-integrability of $K_{ab}(\cdot)$ (see the proof of
Proposition 3.2(i) in \cite{r15}), and integration by parts,
respectively. This further proves (\ref{toprovtheo1}) for $r=s$, hence
also the result.
\end{pf*}
\begin{pf*}{Proof of Theorem \ref{th2}\textup{(ii)} and \textup{(iii)}}
(ii) In view of Theorem \ref{th2}(i), it is sufficient to show that
both asymptotic normality results hold for
$\Delta^*_{\vartheta,f,g;2}$. The result under
$\mathrm{P}^{(n)}_{\vartheta,g}$ then straightforwardly follows from
the multivariate CLT. As for the result under local alternatives
[which, just as the result in part (iii), requires that $g\in\mathcal
{F}_{\mathrm{ulan}}$], it is obtained as usual, by establishing the
joint normality under $\mathrm{P}^{(n)} _{\vartheta,g}$ of
$\log(d\mathrm{P}^{(n)}_{\vartheta+ n^{-1/2}\tau ,g}
/d\mathrm{P}^{(n)}_{\vartheta,g})$ and $\Delta^*_{\vartheta,f,g;2}$,
then applying Le Cam's third lemma; the required joint normality
follows from\vadjust{\goodbreak} a routine application of the classical Cram\'er--Wold
device. (iii) The proof, that is long and tedious, is also a quite
trivial adaptation of the proof of Proposition A.1 in \cite{Ha06c}. We
therefore omit it.
\end{pf*}

\subsection{\texorpdfstring{Proof of Theorem \protect\ref{theotest}}{Proof of Theorem 3.1}}

(i) Applying Corollary \ref{corol}, with
$\check\vartheta:=\break\hat\vartheta_{0}=(\hat\mu',(\operatorname
{vecd}^\circ L_0)')'$ and
$\vartheta:=\vartheta_0=(\mu',(\operatorname{vecd}^\circ L_0)')'$,
entails that\vspace*{-1pt} $\underline{\Delta}^*_{\hat\vartheta_{0},f;2}
=\underline{\Delta}^*_{\vartheta_0,f;2}+o_\mathrm{P}(1)$ as $\ny$ under
$\mathrm{P}^{(n)}_{\vartheta_0,g}$. Consequently, we have that
%
%
\begin{equation} \label{ff}
\underline{Q}_f
=
(\operatorname{vec} \underline{\Delta}^*_{\vartheta_0,f;2})'
(\Gamma^*_{L_0,f;2})^{-1} (\operatorname{vec} \underline{\Delta
}^*_{\vartheta
_0,f;2})+o_\mathrm{P}(1),
\end{equation}
still\vspace*{1pt} as $\ny$, under $\mathrm{P}^{(n)}_{\vartheta_0,g}$---hence also
under $\mathrm{P}^{(n)}_{\vartheta_0+n^{-1/2}\tau,g}$ (from
contiguity). The
result then follows from Theorem \ref{th2}(ii).
%
(ii)
It directly follows from (i) that, under the sequence
of local alternatives $\mathrm{P}^{(n)}_{\vartheta_0+n^{-1/2}\tau,f}$,
$\underline{\phi}_f^{(n)}$ has asymptotic power
$
1
-
\Psi_{p(p-1)}
(
\chi^2_{p(p-1),1-\alpha}
;
\tau_2'
\Gamma^*_{L_0,f;2}
\tau_2
)
$. This establishes the result, since these local powers coincide with
the semiparametrically optimal (at $f$) powers in (\ref{noncentrsemiparam}).


\subsection{\texorpdfstring{Proofs of Lemma \protect\ref{gamrhoestimtheor},
Theorems \protect\ref{theoestim} and \protect\ref{theoestimexplicit}}{Proofs of Lemma 4.1,
Theorems 4.1 and 4.2}}

\mbox{}

\begin{pf*}{Proof of Theorem \ref{theoestim}}
(i)
Fix $\vartheta\in\Theta$ and $g\in\mathcal{F}_{\mathrm{ulan}}$.
From (\ref{definpseudo}), the fact that $\hat\Gamma^*_{\tilde
{L},f;2}-\Gamma^*_{L,f,g;2}=o_\mathrm{P}(1)$
as $\ny$ under $\mathrm{P}^{(n)}_{\vartheta,g}$, and Corollary \ref{corol},
we obtain
%
%
\begin{eqnarray}\label{estimrepresvecdcirc}\quad
\sqrt{n} \operatorname{vecd}^\circ(\hat{\underline{L}}_{f}-L)
&=&
\sqrt{n} \operatorname{vecd}^\circ(\tilde L-L)
+
(\hat\Gamma^*_{\tilde{L},f;2})^{-1}
\underline{\Delta}^*_{\tilde{\vartheta},f;2}
\nonumber
\\
&=&
\sqrt{n} \operatorname{vecd}^\circ(\tilde L-L)
+
(\Gamma^*_{L,f,g;2})^{-1}
\underline{\Delta}^*_{\tilde{\vartheta},f;2}
+o_\mathrm{P}(1)
\nonumber\\
&=&
(\Gamma^*_{L,f,g;2})^{-1}
\underline{\Delta}^*_{\vartheta,f;2}+o_\mathrm{P}(1)
\end{eqnarray}
as $\ny$ under $\mathrm{P}^{(n)}_{\vartheta,g}$. Consequently,
Theorem \ref
{th2}(i) and (ii) entails that, still as $\ny$ under $\mathrm
{P}^{(n)}_{\vartheta,g}$,
%
%
\begin{eqnarray}
\label{represestimcirc}
&&\sqrt{n} \operatorname{vecd}^\circ(\hat{\underline{L}}_{f}-L)\nonumber\\[-8pt]\\[-8pt]
&&\qquad=
(\Gamma^*_{L,f,g;2})^{-1}
{\Delta}^*_{\vartheta,f,g;2}+o_\mathrm{P}(1)
\nonumber\\
\label{estimasymplawvecdcirc}
&&\qquad\stackrel{\mathcal{L}}{\to}
\mathcal{N}_{p(p-1)} ( 0, ( \Gamma^*_{L,f,g;2})^{-1} \Gamma^*_{L,f;2}
(\Gamma^*_{L,f,g;2})^{-1\prime} ).
\end{eqnarray}

Now, by using the fact that $C'(\operatorname{vecd}^\circ
H)=(\operatorname{vec}
H)$ for
any $p\times p$ matrix $H$ with only zero diagonal entries, we have that
$
\sqrt{n} \operatorname{vec}(\hat{\underline{L}}_{f}-L)
=
\sqrt{n} C'\operatorname{vecd}^\circ(\hat{\underline{L}}_{f}-L)
$, so that
(\ref{restim}), (\ref{represestim}) and (\ref{estimasymplawvec})
follow from (\ref{estimrepresvecdcirc}), (\ref{represestimcirc})
and (\ref{estimasymplawvecdcirc}), respectively.


(ii)
The asymptotic covariance matrix of $\sqrt{n} \operatorname{vecd}^\circ
(\hat{\underline{L}}_{f}-L)$,
under $\mathrm{P}^{(n)}_{\vartheta,f}$,
reduces to $(\Gamma^*_{L,f;2})^{-1}$ [let $g=f$ in (\ref
{estimasymplawvecdcirc})], which establishes the result.\vadjust{\goodbreak}
\end{pf*}


To prove Theorem \ref{theoestimexplicit}, we will need the following result.
\begin{Lem}\label{magique}
Fix
$\vartheta=(\mu',(\operatorname{vecd}^\circ L)')'\in\Theta$ and $f,g\in
\mathcal{F}_{\mathrm{ulan}}$. Then
\begin{eqnarray*}
\hspace*{-4pt}&&(I_p\otimes L^{-1})C' (\Gamma^*_{L,f,g;2})^{-1}C(I_p\otimes L^{-1})'
\\[3pt]
\hspace*{-4pt}&&\qquad=
\sum_{r,s=1, r\neq s}^p
\bigl\{
\alpha_{rs}(f,g)
\bigl(e_re_r'\otimes( L_{rs}^2e_r e_r' + e_s e_s' - L_{rs} e_r e_s'
- L_{rs} e_s e_r' )\bigr)
\\[3pt]
\hspace*{-4pt}&&\hspace*{-2pt}\hphantom{\sum_{r,s=1, r\neq s}^p
\bigl\{}
\qquad\quad{} + \beta_{rs}(f,g) \bigl(e_re_s'\otimes( L_{rs}L_{sr} e_r
e_s' - L_{rs} e_r e_r' - L_{sr} e_s e_s' + e_s e_r' )\bigr) \bigr\},
\end{eqnarray*}
where $L_{rs}$ denotes the entry $(r,s)$ of $L$.
\end{Lem}
\begin{pf*}{Proof of Theorem \ref{theoestimexplicit}}
By using again the fact that $C'(\operatorname{vecd}^\circ
H)=(\operatorname{vec} H)$ for any $p\times p$ matrix $H$ with only
zero diagonal entries, and then Lem\-ma~\ref{magique}, we
obtain
\begin{eqnarray*}
&&\operatorname{vec}(\hat{\underline{L}}_{f}-\tilde{L})\\[3pt]
&&\qquad=
C' \operatorname{vecd}^\circ(\hat{\underline{L}}_{f}-\tilde{L})
\\[3pt]
&&\qquad=
\frac{1}{\sqrt{n}}
C' (\hat\Gamma^*_{\tilde{L},f;2})^{-1}
C(I_p\otimes\tilde{L}^{-1})'
\operatorname{vec} \underline{T}_{\tilde{\vartheta},f}
\\[3pt]
&&\qquad=\frac{1}{\sqrt{n}} (I_p\otimes\tilde{L})
\\[3pt]
&&\qquad\quad{} \times \Biggl[ \sum_{r,s=1, r\neq s}^p \bigl\{ \hat\alpha_{rs}(f)
\bigl(e_re_r'\otimes( \tilde{L}_{rs}^2e_r e_r' + e_s e_s' - \tilde {L}_{rs}
e_r e_s' - \tilde{L}_{rs} e_s e_r' )\bigr)
\\[3pt]
&&\hspace*{59pt}\qquad\quad{} + \hat\beta_{rs}(f) \bigl(e_re_s'\otimes(
\tilde{L}_{rs}\tilde{L}_{sr} e_r e_s' - \tilde {L}_{rs} e_r e_r' \\[3pt]
&&\qquad\quad\hspace*{210pt}{}-
\tilde{L}_{sr} e_s e_s' + e_s e_r' )\bigr) \bigr\} \Biggr]\\[3pt]
&&\qquad\quad{}\times\operatorname{vec}
\underline{T}_{\tilde{\vartheta},f}.
\end{eqnarray*}
Since all diagonal entries of $\underline{T}_{\tilde{\vartheta},f}$ are
zeros, we have that
%
%
\begin{eqnarray}\label{ahah}\quad
&&\operatorname{vec}(\hat{\underline{L}}_{f}-\tilde{L})\nonumber\\[3pt]
&&\qquad= \frac{1}{\sqrt{n}} (I_p\otimes\tilde{L})
\nonumber\\[-8pt]\\[-8pt]
&&\qquad\quad{}\times\Biggl[ \sum_{r,s=1, r\neq
s}^p \bigl\{ \hat\alpha_{rs}(f) \bigl(e_re_r'\otimes( e_s e_s' - \tilde{L}_{rs}
e_r e_s' )\bigr)
\nonumber\\
&&\hspace*{59.3pt}\qquad\quad{} + \hat\beta_{rs}(f) \bigl(e_re_s'\otimes( e_s e_r' -
\tilde{L}_{rs} e_r e_r' )\bigr) \bigr\} \Biggr] \operatorname{vec}
\underline{T}_{\tilde{\vartheta},f}.
\nonumber
\end{eqnarray}
%
The identity $(C'\otimes A)(\operatorname{vec}
B)=\operatorname{vec}(ABC)$ then yields
%
%
\[
\operatorname{vec}(\hat{\underline{L}}_{f}-\tilde{L})
=
\frac{1}{\sqrt{n}} (I_p\otimes\tilde{L})
\operatorname{vec}
\Biggl[
\sum_{r,s=1, r\neq s}^p
(\hat N_{f})_{sr} ( e_s e_r' - \tilde{L}_{rs} e_r e_r' ) \Biggr] .
\]
%
%
Hence, we have
%
\begin{eqnarray*}
\hat{\underline{L}}_{f}-\tilde{L} &=& \frac{1}{\sqrt{n}} \tilde{L} \sum_{r,s=1, r\neq
s}^p (\hat N_{f})_{sr} ( e_s e_r' - \tilde{L}_{rs} e_r e_r' )
\\
&=& \frac{1}{\sqrt{n}} \tilde{L} \sum_{r,s=1}^p (\hat N_{f})_{sr} (
e_s e_r' - \tilde{L}_{rs} e_r e_r' ) \\
&=& \frac{1}{\sqrt{n}} \tilde{L} \Biggl(
N_{f} - \sum_{r,s=1}^p \tilde{L}_{rs} (\hat N_{f})_{sr} e_r e_r' \Biggr)
\\
&=& \frac{1}{\sqrt{n}} \tilde{L} \Biggl( \hat N_{f} - \sum_{r=1}^p
(\tilde{L}N_{f})_{rr} e_r e_r' \Biggr) \\
&=& \frac{1}{\sqrt{n}} \tilde{L} \bigl( \hat
N_{f} - \operatorname{diag}(\tilde{L}N_{f}) \bigr),
\end{eqnarray*}
which proves the result.
\end{pf*}
\begin{pf*}{Proof of Lemma \ref{gamrhoestimtheor}}
In this proof, all stochastic convergences are as $\ny$ under $\mathrm
{P}^{(n)} _{\vartheta,g}$. First note that, if $ \check\vartheta:= (
\check\mu', (\operatorname{vecd}^\circ\check{L})')' $ is an arbitrary
locally asymptotically discrete root-$n$ consistent estimator for $
\vartheta= ( \mu', (\operatorname{vecd}^\circ L)')' $, we then have
that
%
%
\begin{eqnarray}
\label{linT}
\operatorname{vec}(\underline{T}_{\check\vartheta,f}-\underline
{T}_{\vartheta,f})
&=&
-
G_{f,g}
(I_p\otimes\check L^{-1})C'
\sqrt{n} \operatorname{vecd}^\circ(\check L-L)\nonumber\\[-8pt]\\[-8pt]
&&{}+o_\mathrm{P}(1)\nonumber
\end{eqnarray}
(compare with Corollary \ref{corol}). Incidentally, note that (\ref
{linT}) implies that
$\operatorname{vec} \underline{T}_{\check\vartheta,f}$ is $O_\mathrm
{P}(1)$ [by
proceeding
exactly as in the proof of Theorem \ref{th2}(i) and (ii), we can indeed
show that, under $\mathrm{P}^{(n)}_{\vartheta,g}$, $\operatorname{vec}
\underline
{T}_{\vartheta,f}$ is asymptotically multinormal, hence stochastically
bounded].

Now, from (\ref{linT}), we obtain
\begin{eqnarray*}
&&\operatorname{vec}(\underline{T}_{\tilde\vartheta_{\lambda}^{\gamma_{rs}},f}-
\underline{T}_{\tilde\vartheta,f})
\\
&&\qquad = - G_{f,g} (I_p\otimes\tilde L^{-1})C' \sqrt{n}
\operatorname{vecd}^\circ(\tilde{L}_{\lambda}^{\gamma_{rs}}-\tilde L)
+o_\mathrm{P}(1)
\\
&&\qquad = - \lambda(\underline{T}_{\tilde\vartheta,f})_{rs} G_{f,g}
(I_p\otimes\tilde L^{-1})C' \operatorname{vecd}^\circ\bigl( \tilde{L}
e_re_s' - \tilde{L} \operatorname{diag}(\tilde{L} e_re_s') \bigr) +o_\mathrm{P}(1),
\end{eqnarray*}
%
which, by using the fact that $C'(\operatorname{vecd}^\circ H)=(
\operatorname{vec} H)$
for any $p\times p$ matrix $H$ with only zero diagonal entries, leads to
\begin{eqnarray*}
&&
\operatorname{vec}(\underline{T}_{\tilde\vartheta_{\lambda}^{\gamma_{rs}},f}-
\underline{T}_{\tilde\vartheta,f})
\\
&&\qquad = - \lambda(\underline{T}_{\tilde\vartheta,f})_{rs} G_{f,g}
(I_p\otimes\tilde L^{-1}) \operatorname{vec}\bigl( \tilde{L} e_re_s' - \tilde{L}
\operatorname{diag}(\tilde{L} e_re_s') \bigr) +o_\mathrm{P}(1)
\\
&&\qquad = - \lambda(\underline{T}_{\tilde\vartheta,f})_{rs} G_{f,g}
\operatorname{vec}\bigl( e_re_s' - \operatorname{diag}(\tilde{L} e_re_s') \bigr)
+o_\mathrm{P}(1).
\end{eqnarray*}
This yields
%
\begin{eqnarray*}
&&
\operatorname{vec}(\underline{T}_{\tilde\vartheta_{\lambda}^{\gamma_{rs}},f}-
\underline{T}_{\tilde\vartheta,f}) \\
&&\qquad= - \lambda
(\underline{T}_{\tilde\vartheta,f})_{rs} G_{f,g} \operatorname{vec}( e_re_s')
+o_\mathrm{P}(1)
\\
&&\qquad= -\lambda(\underline{T}_{\tilde\vartheta,f})_{rs}
\bigl(\gamma_{rs}(f,g) \operatorname{vec}( e_re_s')+\rho_{rs}(f,g)
\operatorname{vec}(
e_se_r')\bigr)\\
&&\qquad\quad{} +o_\mathrm{P}(1) .
\end{eqnarray*}
Premultiplying by $(\underline{T}_{\tilde\vartheta,f})_{rs}
(e_s\otimes
e_r)'$, we then obtain
\[
(\underline{T}_{\tilde\vartheta,f})_{rs}
(\underline{T}_{\tilde\vartheta_{\lambda}^{\gamma_{rs}},f})_{rs}
-
((\underline{T}_{\tilde\vartheta,f})_{rs})^2
=
-\lambda
((\underline{T}_{\tilde\vartheta,f})_{rs})^2
\gamma_{rs}(f,g)
+o_\mathrm{P}(1)
\]
[recall indeed that\vspace*{1pt} $\underline{T}_{\tilde\vartheta,f}=O_\mathrm{P}(1)$],
which establishes the $\gamma$-part of the lemma. The proof of
the $\rho
$-part follows along the exact same lines, but for the fact that the
premultiplication is by $(\underline{T}_{\tilde\vartheta,f})_{sr}
(e_r\otimes e_s)'$.
\end{pf*}
\end{appendix}

\section*{Acknowledgments}

We would like to express our gratitude to the Co-Editor, Professor
Peter B\"{u}hlmann, an Associate Editor and one referee. Their careful
reading of a previous version of the paper and their comments and
suggestions led to a considerable improvement of the present paper.
We are also grateful to Klaus Nordhausen for sending to us the R code
for FastICA authored by Abhijit Mandal.

\begin{supplement}[id=suppA]
\stitle{Further results on tests and a
proof of Theorem \ref{gamrhoestimconsist}\\}
\slink[doi]{10.1214/11-AOS906SUPP} 
\sdatatype{.pdf}
\sfilename{aos906\_supp.pdf}
\sdescription{This supplement provides a
simple explicit expression for the proposed test statistics, derives
local asymptotic powers of the corresponding tests, and presents
simulation results for hypothesis testing. It also gives a proof of
Theorem \ref{gamrhoestimconsist}.}
\end{supplement}

%

\printaddresses

\end{document}